\documentclass[12pt,a4paper]{article}

\usepackage{cite}

\usepackage[a4paper, total={6in, 8.5in}]{geometry} 
\usepackage{hyperref} 
\usepackage{makecell} 

\usepackage{amsmath,amssymb,amsthm}
\usepackage{latexsym}
\usepackage{stmaryrd}
\usepackage{wasysym}
\usepackage{multicol}
\usepackage{multirow}
\usepackage{color}
\usepackage{dirtytalk}
\usepackage{url}
\usepackage{fontawesome} 
\usepackage{csquotes}

\newtheorem{theor}{Theorem}
\newtheorem{claim}[theor]{Claim}

\theoremstyle{definition}

\newtheorem{prop}[theor]{Proposition}
\newtheorem{proposition}[theor]{Proposition}
\newtheorem{cor}[theor]{Corollary}
\newtheorem{deff}{Definition}

\newtheorem*{comment}{Comment}

\newtheorem{discussion}{Discussion}
\theoremstyle{remark}
\newtheorem{rem}{Remark}

\theoremstyle{definition}
\newtheorem{idea}{Idea}
\theoremstyle{definition}

\def\BB{\mathbb}

\usepackage{bm}
\begin{document}

\pagestyle{plain} 

\title{Kontsevich graphs act on Nambu\/--\/Poisson brackets,~IV. 
When the invisible becomes crucial}

\author{Mollie S. Jagoe Brown and  Arthemy V. Kiselev\thanks{Bernoulli Institute for Mathematics, Computer Science and Artificial Intelligence, University of Groningen, P.O. Box 407, 9700 AK Groningen, The Netherlands}%
}

\date{29 September 2025}

\maketitle 


\normalsize 
\begin{abstract}\noindent%
Kontsevich's graphs allow encoding multi\/-\/vectors whose coefficients are differential\/-\/polynomial in the coefficients of a given Poisson bracket on an affine real manifold. 
Encoding formulas by directed graphs adapts to the class of Nambu\/-\/determinant Poisson brackets, yet the graph topology becomes dimension\/-\/specific. 
To inspect whether a given Kontsevich graph cocycle~$\gamma$ acts (non)\/trivially --\,in the second Poisson cohomology\,-- on the space of Nambu brackets, taking a vector field solution $\smash{\vec{X}^\gamma_d}$ from dimension~$d$ does not work in~$d+1$. 
For $2 \leqslant d \leqslant 4$, the action of tetrahedron~$\gamma_3$ on Nambu brackets is known to be a Poisson coboundary, $\smash{\dot{P}} = \llbracket P,\smash{\vec{X}^{\gamma_3}_d} (P)\rrbracket$. We explore which minimal (sub)sets of graphs, encoding (non)\/vanishing objects over $\mathbb{R}^d_{\text{aff}}$, generate the topological data that suffice for a solution $\smash{\vec{X}^{\gamma_3}_{d+1}}$ to appear.
We detect that there can be no solution in higher dimension without invisible graphs that vanish as formulas in $d=3$, but whose descendants do not all vanish over $d=4$.
\end{abstract}


\section{Introduction}\label{introduction}

The Kontsevich graph complex is a differential graded Lie algebra (dgLa) of undirected graphs, equipped with the vertex\/-\/expanding differential 
(see~\cite{ascona,JNMP2017}).
The tetrahedron~$\gamma_3$ 
on 4~vertices and 6~edges 
is an example of a nontrivial graph cocycle; it is `good' to infinitesimally deform Poisson brackets on finite\/-\/dimensional affine real manifolds in such a way that they stay Poisson (see~\cite{ascona} or~\cite{OrMorphism2018,rb
} and references therein).

The $\gamma_3$-flow $\dot{P}=Q^{\gamma_3} 
(P)$ on the space of Poisson structures~$P$ (with $C^3$-differentiable coefficients) 
is a Poisson cocycle in the space of Poisson bi-vectors (with respect to the differential $\partial_P=\llbracket P,\cdot\rrbracket$, here 
$\llbracket\cdot,\cdot\rrbracket$ is the Schouten bracket of multivectors
). Indeed, from the construction of the infinitesimal deformation $Q^{\gamma_3}(P)$, starting with the tetrahedron $\gamma_3$ which has a `good' count of vertices and edges, it is seen that 
\[
\partial_P\bigl(Q^{\gamma_3} 
(P)\bigr)=\llbracket P,Q^{\gamma_3}  
   (P)\rrbracket =0.
\]
The question is 
whether the Poisson cocycle $Q^{\gamma_3} 
(P)\in\ker\partial_P$ is a Poisson coboundary, that is, 
whether the following equation for bi\/-\/vectors,
\[ 
    Q^{\gamma_3}(P)=\llbracket P, \vec{X}^{\gamma_3}(P)\rrbracket,
\] 
is satisfied for some trivialising vector field $\vec{X}^{\gamma_3}(P)$. 
There is \emph{no universal} formula of~$\vec{X}^{\gamma_3}(P)$ that would give us a trivialising vector field solution over all Poisson geometries in all dimensions.

We narrow the set of Poisson structures under study: we consider the class of Nambu\/-\/determinant Poisson brackets on~$\mathbb{R}^d_{\text{aff}}$.
We examine whether the Kontsevich graph flows $\dot{P}=Q^{\gamma_3}_d(P)$ on the spaces of Nambu\/--\/Poisson structures~$P$ over dimensions $2 \leqslant d < 5$ 
are Poisson coboundaries.
So, 
we inspect the existence of 
vector fields $\vec{X}^{\gamma_3}_d(P)$ such that 
\begin{equation}\label{maineq}
  Q^{\gamma_3}_d(P)=\llbracket P, \vec{X}^{\gamma_3}_d(P)\rrbracket.
\end{equation} 
Our overall goal in this class of problems --\,for graph cocycles bigger than the tetrahedron~$\gamma_3$\,-- is to minimize in advance the sets of Kontsevich's directed graphs over which vector field solutions can be constructed (if they exist).
In the papers~\cite{avk,MSJB,fs} we found a formula of $\vec{X}_d^{\gamma_3}(P)$
at~$d=4$. 
Now, we report 
a series of experimental results about the presence or absence of solutions~$\vec{X}_d^{\gamma_3}(P)$ to Eq.~\eqref{maineq} over the (sub)sets of Kontsevich's (micro-)\/graphs chosen over the dimension~$d$ 
in a particular way.\footnote{\label{FootVisIntroGCAOPS}
The immediate sequel V.\ (see~\cite{V}) to I.--III.\ and this paper~IV. is a guide to working with the package \textsf{gcaops}
(\textbf{G}raph \textbf{C}omplex \textbf{A}ction \textbf{O}n \textbf{P}oisson \textbf{S}tructures) for \textsf{SageMath} by Buring \cite{BuringPhD}: 
\url{https://github.com/rburing/gcaops}%
}


\subsection*{Preliminaries}\label{preliminaries}

To keep this paper self-contained, we quote here the necessary preliminaries from~\cite{MSJB}.

\smallskip

Recall that we can express any Poisson bracket in terms of a bi-vector field:
$\{f,g\}=P(f,g)$. We can deform a Poisson bi-vector field~$P$ by a suitable graph cocycle $\gamma$ in the Kontsevich graph complex: 
$\dot{P}=Q^{\gamma}(P)$, where $Q^{\gamma}(P)$ is built of as many copies of~$P$ as there are vertices in $\gamma$, see~\cite{rb}.
The bi\/-\/vector $Q^{\gamma}(P)$ is an infinitesimal symmetry of the Jacobi identity

We let the underlying Poisson manifold be the affine space~$\BB{R}^d$, with Cartesian coordinates given by $\BB{R}^d\ni\boldsymbol{x}$ $=(x_1,x_2,\ldots,x_d)$. 
As for the graph cocycle, we consider the smallest case of the tetrahedron~$\gamma_3$;
to deform, we take the class of Nambu--Poisson brackets on~$\BB{R}^d$.

\begin{deff}[Nambu--Poisson bracket]\label{npdef}
The generalised Nambu-determinant Poisson bracket in dimension $d$ for two smooth functions $f,g\in C^\infty(\BB{R}^d)$ is given by the formula 
\begin{equation}\label{EqNambuPoisson}
\{f,g\}_d(\boldsymbol{x})=\varrho(\boldsymbol{x})\cdot \det \left(
    \dfrac{\partial(f,g,a^1,a^2,\ldots,a^{d-2})}{\partial(x_1,x_2,\ldots,x_d)}
\right)
(\boldsymbol{x}),
\end{equation}
where $a^1,\ldots,a^{d-2}\in C^\infty(\BB{R}^d)$ are Casimirs, which Poisson-commute with any function. The object~$\varrho$ is the inverse density, i.e.\ the coefficient of a $d$-vector field.
\end{deff}

\noindent \textbf{Graph calculus.} 
We solve equation (\ref{maineq}) on the level of formulas. The graph calculus, commonly used in deformation quantisation, is used to obtain these formulas. In essence, graphs represent formulas: the directed edges of the graphs represent derivations which act on the content of vertices. The main convenience of the graph calculus is that formulas change with the dimension, but pictures of graphs do not change. We denote by $\phi$ the map from 
graphs to their formulas obtained by the graph calculus. A detailed explanation of how the graph calculus works can be found in example 1, section 2.2 of \cite{MSJB}. 


\begin{deff}[The sunflower graph]\label{sunflower}
    A linear combination of the above Kontsevich graphs (graphs built of wedges $\smash{\xleftarrow{L}\!\!\bullet\!\!\xrightarrow{R}}$, see~\cite{avk,fs,rb}) can be expressed as the sunflower graph $$\text{sunflower }=\text{ }\raisebox{0pt}[6mm][4mm]{\unitlength=0.4mm
\special{em:linewidth 0.4pt}
\linethickness{0.4pt}
\begin{picture}(17,24)(5,5)
\put(-5,-7){
\begin{picture}(17.00,24.00)
\put(10.00,10.00){\circle*{1}}
\put(17.00,17.0){\circle*{1}}
\put(3.00,17.0){\circle*{1}}
\put(10.00,10.00){\vector(0,-1){7.30}}
\put(17.00,17.00){\vector(-1,0){14.00}}
\put(3.00,17.00){\vector(1,-1){6.67}}
\bezier{30}(3,17)(6.67,13.67)(9.67,10.33)
%
%
\put(17,17){\vector(-1,-1){6.67}}
\bezier{30}(17,17)(13.67,13.67)(10.33,10.33)
\bezier{52}(17.00,17.00)(16.33,23.33)(10.00,24.00)
\bezier{52}(10.00,24.00)(3.67,23.33)(3.00,17.00)
\put(16.8,18.2){\vector(0,-1){1}}
\put(10,17){\oval(18,18)}
\put(10,10){\line(1,0){10}}
\bezier{52}(20,10)(27,10)(21,16)
\put(21,16){\vector(-1,1){0}}
\end{picture}
}\end{picture}}=1\cdot\Gamma_1+2\cdot\Gamma_2.$$ The outer circle means that the outgoing arrow acts on the three vertices via the Leibniz rule. When the arrow acts on the upper two vertices, we obtain two isomorphic graphs, hence the coefficient 2 in the linear combination.
\end{deff}


\begin{prop}[Cf. \cite{ascona}, \cite{Anass2017}]\label{2D}
   There exists a unique (up to 1-dimensional shifts) trivialising vector field in 2D for the $\gamma_3$-flow. It is given by the sunflower graph $$\vec{X}^{\gamma_3}_{d=2}(P)= \phi(\text{sunflower}).$$ The sunflower gives a formula in 2D to solve equation (\ref{maineq}), namely $\dot{P}=Q^{\gamma_3}_{d=2}(P)=\llbracket P,\vec{X}_{d=2}^{\gamma_3}(P)\rrbracket.$ 
\end{prop}

\begin{deff}[Nambu micro-graph]
    Nambu micro-graphs are built using the Nambu--Poisson bracket $P(\varrho,\boldsymbol{a})$ as subgraphs with ordered and directed edges. The vertex of the source of each $P(\varrho,\boldsymbol{a})$ in dimension $d$ contains $\varepsilon^{i_1\hdots i_d}\varrho$, with $d$-many outgoing edges. The first two edges act on the arguments of that bi-vector field subgraph, and the last $d-2$ edges go to the Casimirs $a^1,...,a^{d-2}$.
\end{deff}



\begin{deff}[$d$-descendants] The $d$-descendants of a given $(d'=2)$-dimensional Kontsevich graph is the set of Nambu micro-graphs obtained in the following way. Take a $(d'=2)$-dimensional Kontsevich graph. To each vertex, add $(d-2)$-many Casimirs by $(d-2)$-many outgoing edges. Extend the original incoming arrows to work via the Leibniz rule over the newly added Casimirs. 
\end{deff}

An example of how we can retrieve the $d=3$ descendants of a $d'=2$ Nambu micro-graph can be found in example~4, section~3.2 of~\cite{MSJB}.

\begin{prop}\label{3D} There exists a unique (up to 3-dimensional shifts) trivialising vector field $\vec{X}^{\gamma_3}_{d=3}(P)=\phi(X^{\gamma_3}_{d=3})$ in 3D. It is given as a linear combination over 10 3D-descendants of the 2D sunflower. 
\end{prop} 

\begin{prop}\label{4D} There exists a unique (up to 7-dimensional shifts) trivialising vector field $\vec{X}_{d=4}^{\gamma_3}(P)=\phi(X_{d=4}^{\gamma_3})$ in~4D. 
It is given as a linear combination of 27 \emph{skew pairs} of 1-vector Nambu micro-graphs, where $\text{\emph{skew pair}}=\tfrac{1}{2}\Bigl(\phi\bigl(\Gamma(a^1,a^2)\bigr)-\phi\bigl(\Gamma(a^2,a^1)\bigr)\Bigr)$, and $\Gamma$ belongs to the set of 4D-descendants of the 2D sunflower. 
\footnotesize
\end{prop}

\begin{claim}\label{proj}
    We can project a 4D solution down to a 3D solution by setting the last Casimir equal to the last coordinate: $a^2=w$. Similarly, we can project a 3D solution down to a 2D solution by setting the Casimir equal to the last coordinate: $a^1=z$. We establish that formulas project down to previously found formulas: $$\phi\big(\vec{X}^{\gamma_3}_{d=4}(P)\big)\xrightarrow{a^2=w} \phi\big(\vec{X}^{\gamma_3}_{d=3}(P)\big)\xrightarrow{a^1=z}\phi\big(\vec{X}^{\gamma_3}_{d=2}(P)\big).$$ 
\end{claim}

\begin{claim}\label{nosol3to4}
    There exist solutions in 3D and 4D over the descendants of a 2D solution, the sunflower. But the descendants of known solutions in 3D do not give solutions in~4D.
\end{claim}

\begin{comment}
This would have been practical for reducing computing time. We have 41 3D graphs obtained from expanding the sunflower, and solutions in 3D over only 10 such graphs. Moving up to a higher dimension, it would be ideal to search over descendants of a 3D solution, but we observe this is not possible.
\end{comment}


\section{Preservation and destruction of trivialising vector fields 
in $d=3\mapsto d=4$}

We now examine Claim~\ref{nosol3to4} from~\cite{MSJB}, and we shall expand on the (non)ability to restrict the set of sunflower micro-graphs needed in 3D to obtain a 4D trivialising vector field $\vec{X}^{\gamma_3}_{d=4}(P)$ over their 4D-descendants. 


The most significant obstruction to solving 
problem~\eqref{maineq} in higher dimensions $d\geq5$ is its size. In $d=2$, the plain-text formula of $Q^{\gamma_3}_d(P)$ is a couple of lines long and the solution to \eqref{maineq} can be done by hand. In $d=3$, the plain-text formula of $Q^{\gamma_3}_d(P)$ is a few pages long and a machine computation is required to obtain it and then solve \eqref{maineq}. In $d=4$, the plain-text formula of $Q^{\gamma_3}_d(P)$ is 3GB and a high performance computing cluster is required to obtain it, and solve \eqref{maineq}; see lines~3 and~6 in Table~\ref{table1} 
for precise data on these computational barriers.

\begin{table}[htb]\label{table1}
\centering
\caption{Data regarding the sunflower micro-graphs used in each dimension $d\lesssim 5$ to find a trivialising vector field $\vec{X}^{\gamma_3}_d(P)$.}\label{table1}
\begin{tabular}{|l|l|l|l|l|}
\hline
 & \thead{2D} (see~\cite{ascona,Anass2017}) & \thead{3D} (see~\cite{MSJB}) & \thead{4D} (see~\cite{MSJB}) & \thead{5D} (see App.~\ref{5D}) \\ 
\hline
\thead{Number of sunflower \\ components} & 3 & 48 & 324 & 1280  \\ \hline
\thead{Computation time \\ for one formula} & $\mathcal{O}$(sec) & $\mathcal{O}$(min) & $\mathcal{O}$(min) & 8h \\ \hline
\thead{Number of linearly \\ independent formulas} & 2 & 20 & 123 & ? \\ \hline
\thead{Number of sunflower \\ components in solution} & 2  & 10 & 27 & ? \\ \hline
\thead{Computation time \\ for solving \eqref{maineq}} & $\mathcal{O}$(min) & $\mathcal{O}$(min) & 10h & ? \\ \hline
\end{tabular}
\end{table}


The greatest contribution to solving 
problem~\eqref{maineq} in higher dimensions $d=3,4$ has been to only search for trivialising vector fields $\vec{X}^{\gamma_3}_{d=3,4}(P)$ over $(d=3,4)$-descendants of the 2D sunflower graph (as opposed to searching over all possible graphs in a given dimension; this is explained in sections~3.2 and~3.3 of~\cite{MSJB}). We now attempt to further reduce the number of graphs used in the problem; that is, we attempt to reduce the number of graphs used to obtain formulas which yield a trivialising vector field $\vec{X}^{\gamma_3}_d(P)$, provided it exists at all. Indeed, we still observe huge degeneracy of the number of sunflower micro-graphs needed to express a solution $\vec{X}^{\gamma_3}_d(P)$ for $d=3,4$ in comparison with the total number of sunflower micro-graphs in each dimension. Compare lines~2 and~5 in Table~\ref{table1} for precise data behind this motivation.


Before passing to $d\geq5$, let us examine this reduction of the set of micro-graphs in the $d=3\mapsto d=4$ case. Specifically, we attempt to find a minimal subset of 3D sunflower micro-graphs such that their 4D-descendants yield formulas over which there exists a 4D trivialising vector field $\vec{X}^{\gamma_3}_{d=4}(P)$. Note that we know in advance the space of solutions $\vec{X}^{\gamma_3}_{d=3}(P)$ and $\vec{X}^{\gamma_3}_{d=4}(P)$ over dimensions $d=3,4$, see Propositions 4 and 8 respectively in \cite{MSJB}.

\begin{normalsize}
\begin{idea}[see lines~3,4,5 in Table~\ref{updowntable} on p.~\pageref{updowntable}]\label{idea1}
    Take the 4D-descendants of 3D solutions $\vec{X}_{d=3}^{\gamma_3}(P)$, and search over their formulas for a trivialising vector field $\vec{X}_{d=4}^{\gamma_3}(P)$. 
\end{idea}

\begin{proposition}[code in App.~\ref{A1},\,\ref{A2},\,\ref{A3}]\label{fail1}
    There does \textbf{not} exist a 4D trivialising vector field $\vec{X}^{\gamma_3}_{d=4}(P)$ over the 4D-descendants of the 3D sunflower micro-graph components of known 3D solutions $\vec{X}^{\gamma_3}_{d=3}(P)$. 
\end{proposition} 

\begin{rem}\label{vanishingremark}
    The $(d=4)$-descendants of a given 3D vanishing sunflower micro-graphs are not necessarily vanishing.\footnote{\label{FootEx3Dvanish4Dnot}
For instance, the 3D vanishing sunflower micro-graph given by the encoding (0,1,4;1,6,5;4,5,6) has a 4D non-vanishing sunflower micro-graph descendant given by the encoding (0,1,4,7;1,6,5,8;4,8,6,9); see section~2.2 in~\cite{MSJB} for an explanation of the graph encodings.
}
For our future use, here are the encodings of all 3D vanishing sunflower micro-graphs; the first two are isomorphic. (By definition, a graph is \textit{zero} when it has a symmetry under which it equals minus itself, hence its formula is identically zero.)
{    \small
\begin{multicols}{3}
\begin{itemize}
\item[1.](0,1,4;1,6,5;4,5,6)
\item[2.](0,1,4;4,6,5;1,5,6)
\item[3.](0,2,4;1,3,5;4,2,6)
\item[4.](0,5,4;4,3,5;1,2,6)
\item[5.](0,5,4;1,3,5;4,2,6) 
\item[6.](0,5,4;4,6,5;1,2,6)
\item[7.](0,5,4;1,6,5;4,2,6)
\item[8.](0,2,4;1,3,5;4,5,6)
\item[9.](0,2,4;4,6,5;1,5,6)
\item[10.](0,2,4;1,6,5;4,5,6)
\item[11.](0,5,4;1,3,5;4,5,6)
\item[12.](0,5,4;4,3,5;4,5,6) \textit{zero}
\item[13.](0,5,4;4,6,5;4,5,6) \textit{zero}
\end{itemize}
\end{multicols}

}
The encodings of all 4D-descendants of 3D vanishing sunflower micro-graphs are given in  Appendix~\ref{A8}. 
\end{rem}

\begin{idea}\label{idea2}
    Take the 4D-descendants of the 3D vanishing sunflower micro-graphs; adjoin these to the set of 4D-descendants of the 3D sunflower micro-graph components of known 3D solutions $\vec{X}^{\gamma_3}_{d=3}(P)$, and search over their formulas for a trivialising vector field $\vec{X}_{d=4}^{\gamma_3}(P)$.
\end{idea}

\begin{proposition}[code in App~\ref{A1},\,\ref{A2},\,\ref{A3}]\label{fail2}
    There does \textbf{not} exist a 4D trivialising vector field $\vec{X}^{\gamma_3}_{d=4}$ over the union of the 4D-descendants of the 3D vanishing sunflower micro-graphs and the 4D-descendants of the 3D sunflower micro-graph components of known 3D solutions $\vec{X}^{\gamma_3}_{d=3}(P)$. 
\end{proposition}

We still lose valuable micro-graphs in the 3D$\mapsto$4D step. 
Let us recall that by Claim~\ref{proj}, the formula of a 4D trivialising vector field $\vec{X}^{\gamma_3}_{d=4}(P)$ projects down to the formula of a 3D trivialising vector field by setting the last Casimir in 4D equal to the last coordinate in 4D, $a^2=w$. 

\begin{idea}\label{4dto3d}
    Project a known 4D trivialising vector field $\vec{X}^{\gamma_3}_{d=4}(P)$ down to a 3D trivialising vector field $\vec{X}_{d=3}^{\gamma_3}(P)$ under $a^2=w$; identify the 3D sunflower micro-graph components of this $\vec{X}^{\gamma_3}_{d=3}(P)$, and search for a trivialising vector field $\vec{X}^{\gamma_3}_{d=4}(P)$ over the formulas of their 4D-descendants.
\end{idea}

\begin{proposition}[code in App.~\ref{A4}]\label{fail3}
    There does \textbf{not} exist a 4D trivialising vector field $\vec{X}^{\gamma_3}_{d=4}(P)$ over the 4D-descendants of the 3D sunflower micro-graph components of the 3D solution $\vec{X}^{\gamma_3}_{d=3}(P)$ obtained by projecting a known 4D solution $\vec{X}^{\gamma_3}_{d=4}(P)$ down under $a^2=w$. 
\end{proposition}

The set of 3D sunflower micro-graphs which we take is still not sufficient to give enough 4D-descendants as to yield a 4D trivialising vector field $\vec{X}_{d=4}^{\gamma_3}(P)$. So, we now turn to the twenty 3D sunflower micro-graphs which yield linearly independent formulas. Their encodings are: 
{
\small
\begin{multicols}{4}
\begin{itemize}
\item[1.](0,1,4;1,3,5;1,2,6)  
\item[2.](0,1,4;1,6,5;1,2,6)
\item[3.](0,1,4;1,3,5;4,2,6)
\item[4.](0,1,4;1,6,5;4,2,6)
\item[5.](0,1,4;4,6,5;1,2,6)
\item[6.](0,1,4;4,3,5;4,2,6)
\item[7.](0,1,4;4,6,5;4,2,6)
\item[8.](0,1,4;1,6,5;1,5,6)
\item[9.](0,1,4;4,6,5;4,5,6)
\item[10.](0,2,4;1,3,5;1,2,6)
\item[11.](0,2,4;1,6,5;1,2,6)
\item[12.](0,2,4;4,3,5;4,2,6)
\item[13.](0,2,4;4,6,5;4,2,6)
\item[14.](0,2,4;1,3,5;1,5,6)
\item[15.](0,2,4;1,6,5;1,5,6)
\item[16.](0,2,4;4,3,5;1,5,6)
\item[17.](0,5,4;1,3,5;1,2,6)
\item[18.](0,5,4;4,3,5;4,2,6)
\item[19.](0,5,4;1,3,5;1,5,6)
\item[20.](0,5,4;1,6,5;1,5,6)
\end{itemize}
\end{multicols}

}
\normalsize

\begin{rem}\label{syn}
  There exist \emph{synonyms}: different graphs can give the same formulas. The twenty 3D sunflower micro-graphs yielding linearly independent formulas are chosen by the internal algorithms of \textsf{SageMath}. Some \emph{other} choice of 20 basic graphs yielding linearly independent formulas could alter our conclusions in what follows.   
\end{rem}

\begin{idea}\label{ind}
    Take the 4D-descendants of the twenty 3D sunflower micro-graphs which yield linearly independent formulas, and search over their formulas for a 4D trivialising vector field $\vec{X}_{d=4}^{\gamma_3}(P)$. 
\end{idea}

\begin{proposition}[code in App.~\ref{A5}]\label{fail4}
    There does \textbf{not} exist a 4D trivialising vector field $\vec{X}_{d=4}^{\gamma_3}(P)$ over the 4D-descendants of the 3D sunflower micro-graphs which yield linearly independent formulas. 
\end{proposition}

Clearly, there must be a source of relevant micro-graphs which we have been missing so far.

\begin{idea}\label{manual}
    Draw the 4D sunflower micro-graphs used in skew-symmetrised formulas (with respect to $a^1\rightleftarrows a^2$) giving a 4D trivialising vector field $\vec{X}^{\gamma_3}_{d=4}(P)$. Compare these drawings with the full set of the 3D sunflower micro-graphs: identify which 3D sunflower micro-graphs we need to obtain the 4D sunflower micro-graphs that actually occur in the skew pairs in the 4D trivialising vector field $\vec{X}^{\gamma_3}_{d=4}(P)$. 
\end{idea}

\begin{proposition}\label{17_3D}
    We find that seventeen 
3D sunflower micro-graphs are needed to obtain a sufficiently large set of 4D-descendants so as to include the 4D sunflower micro-graphs which yield a known 4D trivialising vector field $\vec{X}_{d=4}^{\gamma_3}(P)$. Of these 17 micro-graphs over 3D, we have 3 vanishing micro-graphs (numbered 7,10,13 in Remark~\ref{vanishingremark}); 12 belonging to the \textsf{SageMath}-chosen basis of 20 giving linearly independent formulas (numbered 1,4,6,7,8,9,10,14,15,16,17,20 above Remark~\ref{syn}), and 2 remain (they are non-vanishing and not in the basis: (0,5,4;1,6,5;1,2,6) and (0,5,4;4,6,5;4,2,6)). That is, these two micro-graphs are non-vanishing and have formulas obtained by linear combinations of the 20 linearly independent formulas. 
\end{proposition}

\begin{discussion}\label{graphlanguageresists} The two non-vanishing graphs whose formulas are obtained as linear combinations of the 20 linearly independent formulas offer interesting insight into the resistance of the graph calculus from one dimension to the next, in this case from 3D to 4D.
We pose the question: does the graph calculus preserve linear combinations from dimension~$d$ to $d+1$? 
Assume that the formula of a non-vanishing 3D micro-graph $\Gamma_{3D}$ is expressed as a linear combination of linearly independent formulas obtained from 3D micro-graphs $\Gamma^i_{3D}$: 
\[
\phi(\Gamma_{3D})=\sum\alpha_i\cdot\phi(\Gamma^i_{3D}),
\]
with scalars $\alpha_i$. Are the formulas of the 4D-descendants of $\Gamma_{3D}$, denoted by $\Gamma_{3D\mapsto4D}$, expressed as linear combinations of the formulas of the 4D-descendants $\Gamma_{3D\to4D}^i$ of the 3D micro-graphs $\Gamma^i_{3D}$ which yield linearly independent formulas? That is, is it always true that 
for some scalars $\beta_i$,  
\begin{equation}\label{EqLinDepDLinDepHigher}
\phi(\Gamma_{3D})=\sum\alpha_i\cdot\phi(\Gamma^i_{3D}) \quad \overset{?}{\Longrightarrow} \quad \phi(\Gamma_{3D\mapsto 4D})=\sum\beta_i\cdot\phi(\Gamma^i_{3D\to 4D})\,?
\end{equation}
\end{discussion}

\begin{idea}\label{success_idea}
    Take the 4D-descendants of the seventeen 3D sunflower micro-graphs described in proposition~\ref{17_3D} which, we know, contain the 4D sunflower micro-graphs that yield a 4D trivialising vector field $\vec{X}_{d=4}^{\gamma_3}(P)$. Now, search over their formulas for a 4D trivialising vector field $\vec{X}_{d=4}^{\gamma_3}(P)$. 
\end{idea}

\begin{proposition}[code in App.~\ref{A6}]\label{17_3D_give_4D_sol}
    There \textbf{does exist} a 4D trivialising vector field $\vec{X}^{\gamma_3}_{d=4}(P)$ over the 4D-descendants of the seventeen 3D sunflower micro-graphs needed to obtain a sufficiently large set of 4D-descendants so as to include the 4D sunflower micro-graphs which yield a known 4D trivialising vector field $\vec{X}^{\gamma_3}_{d=4}(P)$\,!
\end{proposition}

\begin{discussion}\label{discussion} Intuitively, this result makes sense. We retrieved the 3D graphs needed to obtain the 4D graphs which themselves are necessary 
for a 4D trivialising vector field, so naturally we should be able to obtain the desired 4D solution. 

$\bullet$ The strength of this result can be compared with Idea~\ref{4dto3d}, in which we projected the formula of $\vec{X}^{\gamma_3}_{d=4}(P)$ down to a $\vec{X}^{\gamma_3}_{d=3}(P)$, then searched over the 4D-descendants of this 3D vector field. Critical information can get lost in the graph-to-formula step. 

$\bullet$ This result offers us an answer to the question posed in Discussion~\ref{graphlanguageresists}: yes, we observe that the graph calculus resists the dimensional shift 
3D${}\mapsto{}$4D. That is, the two 3D sunflower micro-graphs whose formulas are linear combinations of linearly independent formulas have 4D-descendants whose formulas are \emph{also} linear combinations of linearly independent formulas.\footnote{\label{FootFirstFormulasThenLinInd}
We know this because in the code, we first compute all formulas of the graphs, then identify which are linearly independent, then search for a solution over these linearly independent formulas.}
This 
is 
remarkable because we observe the 
preservation despite the explosion of the number and nature of 4D-descendants due to the Leibniz rule expansion with respect to the new Casimir. 

$\bullet$ 
The drawback of this result is that we used the \emph{known} 4D trivialising vector field $\vec{X}^{\gamma_3}_{d=4}(P)$ and \emph{graphically} reversed the dimensional step 3D$\mapsto$4D in order to find a restricted set of 3D sunflower micro-graphs. So, this result rather offers insight into the calculus 
of the graph at hand, instead of offering us a new method for restricting the set of $d$-dimensional sunflower micro-graphs in the dimensional shift $d\mapsto d+1$.
\end{discussion}

\begin{idea}\label{general}
    Take the 210 (rather than all of 324) 4D-descendants of the full set of 20 basic 3D sunflower micro-graphs and the 4D-descendants of the full set of vanishing 3D sunflower micro-graphs. Search over the formulas of these 4D-descendants for a 4D trivialising vector field $\vec{X}_{d=4}^{\gamma_3}(P)$.
\end{idea}

\begin{proposition}[code in App.~\ref{A7}]\label{lin_ind_and_vanishing_to_4d}
    There \textbf{exists} a 4D trivialising vector field $\vec{X}_{d=4}^{\gamma_3}(P)$ over the formulas of the 4D-descendants of the 20 basic \textsf{SageMath}-chosen 3D sunflower micro-graphs and the 4D-descendants of the 3D vanishing micro-graphs.
\end{proposition}


\begin{cor}\label{3simps}
    The three simplifications $\#1,\#2,\#3$ outlined in~\cite{MSJB} reduced the size of our 4D problem \eqref{maineq} approximately 300 times: 
\[
\infty\xrightarrow{\text{graphs}}19\text{ }957 \xrightarrow{\#1,\#2} 324 \xrightarrow{\#3} 64.
\]
Here, $\infty$ is the number of all possible formulas; 19,957 is the number of all 1-vector micro-graphs built of 3 Levi--Civita symbols, 3 Casimirs $a^1$ and 3 Casimirs $a^2$; 324 is the number of 4D-descendants of the 2D sunflower; 64 is the number of skew pairs obtained from the 123 linearly independent formulas of the 324 4D-descendants. Now, we reduce the size of the problem even further: instead of looking over 324 4D-descendants of the 2D sunflower, we only need to look over 210 4D-descendants of the 3D sunflower micro-graphs whose formulas are linearly independent, and the 3D vanishing sunflower micro-graphs. 
\end{cor}

The summary of our results 
is contained in Table~\ref{updowntable} on p.~\pageref{updowntable}. 
\begin{table}[h!]
\begin{small}
\begin{center}
\caption{(Non)success of narrowing 
the set of 3D sunflower micro-graphs 
whose 4D-descendants are taken
to obtain a 4D trivialising vector field~$\vec{X}^{\gamma_3}_{d=4}(P)$.}\label{updowntable}
\begin{tabular}{|l|l|l|l|l|}
\hline
 & \thead{Number of 3D \\ sunflower graphs} & \thead{Number of \\4D-descendants} & \thead{Number of lin.
 \\independent\\ formulas} & \thead{Solution in 4D?} \\ \hline
\thead{\textbf{Full 3D sunflower}~\cite{MSJB}} & \textbf{48${}^{\text{*}}$} & \textbf{324} & \textbf{123} & \textbf{yes!}\quad \cite{MSJB}  \\ \hline
\thead{3D solution $\#$1} & 9 & 42 &  & no\quad (\ref{A1}) \\ \hline
\thead{3D solution $\#$2} & 10 & 46 &  & no\quad (\ref{A2}) \\ \hline
\thead{3D solution $\#$3} & 12 & 58 &  & no\quad (\ref{A3}) \\ \hline
\thead{3D vanishing} & 13 & 118${}^{\text{**}}$ &  & NA${}^{\text{***}}$ (\ref{A8}) \\ \hline
\thead{3D solution $\#$1 \\+ 3D vanishing} & 22 & 160 & 72 & no\quad (\ref{A1}) \\ \hline
\thead{3D solution $\#$2 \\+ 3D vanishing} & 23 & 164 & 75 & no\quad (\ref{A2}) \\ \hline
\thead{3D solution $\#$3 \\+ 3D vanishing} & 25 & 176 & 86 & no\quad (\ref{A3}) \\ \hline
\thead{3D solution projected \\ from 4D solution} & 12 &  &  & no\quad (\ref{A4}) \\ \hline
\thead{3D solution projected \\from 4D solution \\+ 3D vanishing} & 25 & 176 & 86 & no \quad(\ref{A4}) \\ \hline
\thead{3D with linearly \\independent formulas} & 20 & 92 & 81
& no\quad (\ref{A5}) \\ \hline
\thead{3D micro-graphs which give\\  the 4D-descendants in \\the 4D solution} & 17 & 110 & 76 & \textbf{yes!}\quad (\ref{A6})\\ \hline
\thead{3D with linearly independent \\formulas + 3D vanishing} & 33 & 210 & 112 & \textbf{yes!}\quad (\ref{A7}) \\ \hline
\end{tabular}\\[3pt]
\end{center}

\end{small}
\noindent\rule{0.8in}{0.7pt}\\
${}^{\text{*}}$\quad Only 41 of the 3D sunflower micro-graphs are non-isomorphic.\\[0.5pt]
${}^{\text{**}}$\quad Note how large this number is compared with the number of 4D-descendants in the three rows above. %
\\[0.5pt]
${}^{\text{***}}$ There cannot exist a 4D solution over this set because the graphs involved do not offer a high enough degree of differentiation of the coefficient~$\varrho$ (in Nambu\/--\/Poisson brackets~\eqref{EqNambuPoisson}) to satisfy~\eqref{maineq}; namely, the embeddings (see Definition~3 in~\cite{fs}) of the 2D sunflower into 4D must appear in any 4D solution.
\end{table}
Each line of the table 
reports on our attempt to restrict the set of 3D sunflower micro-graphs\footnote{\label{FootWeTookSunflower}
As we pointed out in~\cite{avk}, our ability to find the trivialising vector field~$\vec{X}^{\gamma_3}_{d=4}(P)$ solution of~\eqref{maineq} over dimension~$d=4$ was based on an 
irrational idea $(i)$ to only search for 
$\vec{X}^{\gamma_3}_{d=4}(P)$ over the descendants of the sunflower graph.
Now, we detect that to obtain this solution in~4D, it suffices to know a small subset of micro\/-\/graphs for the solution in~3D, and $(ii)$ adjoin the set of 3D vanishing sunflower micro-graphs.%
}
such that their 4D-descendants would yield a 4D trivialising vector field~$\vec{X}^{\gamma_3}_{d=4}(P)$. We ultimately found two new (i.e.\ smaller than taking all) such sets: they are described in Proposition~\ref{17_3D_give_4D_sol} and Proposition~\ref{lin_ind_and_vanishing_to_4d}.
As a by\/-\/product, we observe 
that 
linear dependence 
of formulas obtained from graphs is preserved under the dimensional step $d\mapsto d+1$:
expanding --\,from 3D to 4D, via the Leibniz rule for each arrow that acted on the Casimir in 3D\,-- the 3D micro\/-\/graphs which gave linearly dependent formulas, and translating the linear combination of 4D micro\/-\/graphs into formulas, we assert their linear dependence, 
see~\ref{EqLinDepDLinDepHigher} on p.~\pageref{EqLinDepDLinDepHigher}. 

\subsubsection*{Acknowledgements}
The authors thank the organizers of the ISQS29 conference on integrable systems and quantum symmetries (7--11 July 2025 in CVUT Prague, CZ) for a dynamic and welcoming atmosphere. The first author thanks the Center for Information Technology of the University of Groningen for access to the High Performance Computing cluster, H\'abr\'ok. 
The authors thank F.\,
Schipper and R.\,
Buring for fruitful 
discussions and advice, as well as for their help with coding in \textsf{gcaops} and optimising the code.
The participation of M.S.~Jagoe Brown in the ISQS29 was
supported by the Master's Research Project funds at the Bernoulli Institute, University of Groningen; that of A.V.~Kiselev was supported by project~135110.

\newpage
\appendix
\section{Program codes in \textsf{SageMath} / \textsf{gcaops} and their output}\label{appendix}
Here, the order of the scripts follows the order of ideas and propositions in this paper. The respective proposition of each result is given.

\subsection{Solution n.1 over sunflower micro-graphs in 3D does not extend to a 4D solution}\label{A1}

This result is in Ideas \ref{idea1},\ref{idea2}, Propositions \ref{fail1},\ref{fail2}, and lines 3,4,5 and 7,8,9 in Table~\ref{updowntable} on p.~\pageref{updowntable}. This is because in this script, we include the 4D-descendants of the vanishing 3D sunflower micro-graphs directly.

The script:
\tiny
\begin{verbatim}
NPROCS = 32

# Take the #9 over 3D, plus the #13 vanishing (out of full #48 sunflower), extend to 4D: is there a solution?

import warnings
warnings.filterwarnings("ignore", category=DeprecationWarning)

# Graphs for the vector field $X$

X_graph_encodings = [(0,1,4,7,1,3,5,8,1,2,6,9),(0,1,4,7,4,3,5,8,4,2,6,9),(0,1,4,7,7,3,5,8,4,2,6,9),(0,1,4,7,7,3,5,8,7,2,6,9),
(0,1,4,7,4,3,5,8,7,2,6,9),(0,1,4,7,4,6,5,8,4,2,6,9),(0,1,4,7,7,6,5,8,4,2,6,9),(0,1,4,7,7,6,5,8,7,2,6,9),(0,1,4,7,4,6,5,8,7,2,6,9),
(0,1,4,7,4,9,5,8,7,2,6,9),(0,1,4,7,7,9,5,8,7,2,6,9),(0,1,4,7,7,9,5,8,4,2,6,9),(0,1,4,7,4,9,5,8,4,2,6,9),
(0,5,4,7,1,3,5,8,1,2,6,9),(0,8,4,7,1,3,5,8,1,2,6,9),(0,5,4,7,1,6,5,8,1,2,6,9),(0,8,4,7,1,6,5,8,1,2,6,9),(0,5,4,7,1,9,5,8,1,2,6,9),
(0,8,4,7,1,9,5,8,1,2,6,9),(0,2,4,7,1,3,5,8,1,2,6,9),(0,2,4,7,1,3,5,8,1,5,6,9),(0,2,4,7,1,3,5,8,1,8,6,9),
(0,2,4,7,4,3,5,8,1,5,6,9),(0,2,4,7,7,3,5,8,1,5,6,9),(0,2,4,7,7,3,5,8,1,8,6,9),(0,2,4,7,4,3,5,8,1,8,6,9)]
for i in [4,7]:
   for j in [6,9]:
       for k in [4,7]:
           for l in [5,8]:
               X_graph_encodings.append((0,1,4,7,i,j,5,8,k,l,6,9))

vanishing_expanded=[(0,2,4,7,1,3,5,8,4,2,6,9),(0,2,4,7,1,3,5,8,7,2,6,9),(0,5,4,7,4,3,5,8,1,2,6,9),(0,8,4,7,4,3,5,8,1,2,6,9),
(0,8,4,7,7,3,5,8,1,2,6,9),(0,5,4,7,7,3,5,8,1,2,6,9),(0,5,4,7,1,3,5,8,4,2,6,9),(0,8,4,7,1,3,5,8,4,2,6,9),
(0,8,4,7,1,3,5,8,7,2,6,9),(0,5,4,7,1,3,5,8,7,2,6,9),(0,2,4,7,1,3,5,8,4,5,6,9),(0,2,4,7,1,3,5,8,7,5,6,9),
(0,2,4,7,1,3,5,8,7,8,6,9),(0,2,4,7,1,3,5,8,4,8,6,9)]

for i in [4,7]:
    for j in [5,8]:
        for k in [6,9]:
            vanishing_expanded.append((0,1,4,7,1,k,5,8,i,j,6,9))
            vanishing_expanded.append((0,2,4,7,1,k,5,8,i,j,6,9))
            vanishing_expanded.append((0,2,4,7,i,k,5,8,1,j,6,9))
            vanishing_expanded.append((0,j,4,7,i,k,5,8,1,2,6,9))
            vanishing_expanded.append((0,j,4,7,1,k,5,8,i,2,6,9))
            vanishing_expanded.append((0,1,4,7,i,k,5,8,1,j,6,9))
            
for i in [4,7]:
    for j1 in [5,8]:
        for j2 in [5,8]:
            vanishing_expanded.append((0,j1,4,7,1,3,5,8,i,j2,6,9))
            
for i1 in [4,7]:
    for i2 in [4,7]:
        for j1 in [5,8]:
            for j2 in [5,8]:
                vanishing_expanded.append((0,j1,4,7,i1,3,5,8,i2,j2,6,9))
                
for i1 in [4,7]:
    for i2 in [4,7]:
        for j1 in [5,8]:
            for j2 in [5,8]:
                for k in [6,9]:
                    vanishing_expanded.append((0,j1,4,7,i1,k,5,8,i2,j2,6,9))

X_graph_encodings += vanishing_expanded


# Convert the encodings to graphs:

from gcaops.graph.formality_graph import FormalityGraph
def encoding_to_graph(encoding):
    targets = [encoding[0:4], encoding[4:8], encoding[8:12]]
    edges = sum([[(k+1,v) for v in t] for (k,t) in enumerate(targets)], [])
    return FormalityGraph(1, 9, edges)

X_graphs = [encoding_to_graph(e) for e in X_graph_encodings]
print("Number of graphs in X:", len(X_graphs), flush=True)

# Formulas of graphs for $X$

# We want to evaluate each graph to a formula, inserting three copies of $\varrho\varepsilon^{ijk}$ 
# and three copies of the Casimir $a$ into the aerial vertices.

# First define the even coordinates $x,y,z$ (base), $\varrho, a$ (fibre) and odd coordinates $\xi_0, \xi_1, \xi_2$:

from gcaops.algebra.differential_polynomial_ring import DifferentialPolynomialRing
D4 = DifferentialPolynomialRing(QQ, ('rho','a1','a2'), ('x','y','z','w'), max_differential_orders=[3+1,1+3+1,1+3+1])
rho, a1, a2 = D4.fibre_variables()
x,y,z,w = D4.base_variables()
even_coords = [x,y,z,w]

from gcaops.algebra.superfunction_algebra import SuperfunctionAlgebra

S4.<xi0,xi1,xi2,xi3> = SuperfunctionAlgebra(D4, D4.base_variables())
xi = S4.gens()
odd_coords = xi
epsilon = xi[0]*xi[1]*xi[2]*xi[3] # Levi-Civita tensor

# Now evaluate each graph to a formula:

import itertools
from multiprocessing import Pool

def evaluate_graph(g):
    E = x*xi[0] + y*xi[1] + z*xi[2] + w*xi[3] # Euler vector field, to insert into ground vertex. Incoming derivative d/dx^i will result in xi[i].
    result = S4.zero()
    for index_choice in itertools.product(itertools.permutations(range(4)), repeat=3):
        sign = epsilon[index_choice[0]] * epsilon[index_choice[1]] * epsilon[index_choice[2]]
        # NOTE: This assumes the ground vertex is labeled 0, the vertices of out-degree 4 are labeled 1, 2, 3, 
        # and the Casimirs are labeled 4, 5, 6; 7, 8, 9.
        vertex_content = [E, S4(rho), S4(rho), S4(rho), S4(a1), S4(a1), S4(a1), S4(a2), S4(a2), S4(a2)]
        for ((source, target), index) in zip(g.edges(), sum(map(list, index_choice), [])):
            vertex_content[target] = vertex_content[target].derivative(even_coords[index])
        result += sign * prod(vertex_content)
    return result

print("Calculating X_formulas", flush=True)
X_formulas = []
with Pool(processes=NPROCS) as pool:
    X_formulas = list(pool.imap(evaluate_graph, X_graphs))
# NOTE: Necessary fixup of parents after multiprocessing.
X_formulas = [S4(X_formula) for X_formula in X_formulas]
print("Calculated X_formulas", flush=True)

# Each resulting formula is of the form $X_g = \sum_{i=1}^3 X_g^i \xi_i$ where $X_g^i$ is a differential polynomial depending on the graph $g$.

# Relations between (formulas of) graphs in $X$

# For each $i = 1,2,3,4$ we collect the distinct differential monomials in $X_g^i$ ranging over all $g$:

X_monomial_basis = [set([]) for i in range(4)]
for i in range(4):
    for X in X_formulas:
        X_monomial_basis[i] |= set(X[i].monomials())
X_monomial_basis = [list(b) for b in X_monomial_basis]
X_monomial_index = [{m : k for k, m in enumerate(b)} for b in X_monomial_basis]

print("Number of monomials in components of X:", [len(b) for b in X_monomial_basis], flush=True)

X_monomial_count = sum(len(b) for b in X_monomial_basis)

# Now the graph-to-formula evaluation can be expressed as a matrix:

print("Calculating X_evaluation_matrix", flush=True)
X_evaluation_matrix = matrix(QQ, X_monomial_count, len(X_graphs), sparse=True)
for i in range(len(X_graphs)):
    v = vector(QQ, X_monomial_count, sparse=True)
    index_shift = 0
    for j in range(4):
        f = X_formulas[i][j]
        for coeff, monomial in zip(f.coefficients(), f.monomials()):
            monomial_index = X_monomial_index[j][monomial]
            v[index_shift + monomial_index] = coeff
        index_shift += len(X_monomial_basis[j])
    X_evaluation_matrix.set_column(i, v)
print("Calculated X_evaluation_matrix", flush=True)

# Here is the number of linearly independent graphs:

print("Number of linearly independent graphs in X:", X_evaluation_matrix.rank(), flush=True)

# We collect the linearly independent graphs and their formulas in new shorter lists:

pivots = X_evaluation_matrix.pivots()
print("Maximal subset of linearly independent graphs in X:", list(pivots), flush=True)
X_graphs_independent = [X_graphs[k] for k in pivots]
X_formulas_independent = [X_formulas[k] for k in pivots]

# The tetrahedron operation $\operatorname{Op}(\gamma_3)$

from gcaops.graph.undirected_graph_complex import UndirectedGraphComplex

GC = UndirectedGraphComplex(QQ, implementation='vector', sparse=True)
tetrahedron = GC.cohomology_basis(4, 6)[0]

from gcaops.graph.directed_graph_complex import DirectedGraphComplex

dGC = DirectedGraphComplex(QQ, implementation='vector')
tetrahedron_oriented = dGC(tetrahedron)
tetrahedron_oriented_filtered = tetrahedron_oriented.filter(max_out_degree=2)
tetrahedron_operation = S4.graph_operation(tetrahedron_oriented_filtered)

# The flow $\dot{a}_i = 4\cdot\operatorname{Op}(\gamma_3)(P,P,P,a_i)$

# Let $P$ be the 4D Nambu–Poisson bracket with pre-factor $\varrho$ and Casimir functions $a_1, a_2$:

P = -(rho*epsilon).bracket(a1).bracket(a2)

# We calculate $\dot{a}_i = 4\cdot\operatorname{Op}(\gamma_3)(P,P,P,a_i)$:

def casimir_flow(f):
    return 4*tetrahedron_operation(P,P,P,f)

a = [S4(a1), S4(a2)]

print("Calculating adot", flush=True)
adot = [casimir_flow(a_i) for a_i in a]
# adot = []
# with Pool(processes=2) as pool:
#    adot = list(pool.imap(casimir_flow, a))
# NOTE: Necessary fixup of parents after multiprocessing.
adot = [S4(D4(adot_i[()])) for adot_i in adot]
print("Calculated adot", flush=True)

# We calculate $X_g(a_i) = [\![X_g, a_i]\!]$ for $i=1,2$:

# TODO: Parallelize
print("Calculating X_a_formulas", flush=True)
X_a_formulas = [[X_formula.bracket(f) for X_formula in X_formulas_independent] for f in a]
print("Calculated X_a_formulas", flush=True)

# We express $g \mapsto X_g(a_i)$ as a sparse matrix for $i=1,2$:

X_a_basis = [set(f[()].monomials()) for f in adot]
for k in range(len(a)):
    for X_a_formula in X_a_formulas[k]:
        X_a_basis[k] |= set(X_a_formula[()].monomials())
X_a_basis = [list(B) for B in X_a_basis]
print("Number of monomials in X_a_basis:", len(X_a_basis[0]), len(X_a_basis[1]), flush=True)

print("Calculating X_a_evaluation_matrix", flush=True)
X_a_monomial_index = [{m : k for k, m in enumerate(B)} for B in X_a_basis]
X_a_evaluation_matrix = [matrix(QQ, len(B), len(X_graphs_independent), sparse=True) for B in X_a_basis]
for i in range(len(a)):
    for j in range(len(X_graphs_independent)):
        v = vector(QQ, len(X_a_basis[i]), sparse=True)
        f = X_a_formulas[i][j][()]
        for coeff, monomial in zip(f.coefficients(), f.monomials()):
            monomial_index = X_a_monomial_index[i][monomial]
            v[monomial_index] = coeff
        X_a_evaluation_matrix[i].set_column(j, v)
print("Calculated X_a_evaluation_matrix", flush=True)

# We express $\dot{a}_i$ as a vector for $i=1,2$:

print("Calculating adot_vector", flush=True)
adot_vector = [vector(QQ, len(B)) for B in X_a_basis]
for i in range(len(a)):
    f = adot[i][()]
    for coeff, monomial in zip(f.coefficients(), f.monomials()):
        monomial_index = X_a_monomial_index[i][monomial]
        adot_vector[i][monomial_index] = coeff
print("Calculated adot_vector", flush=True)

# The flow $\dot{\varrho}$

print("Calculating Q_tetra", flush=True)
Q_tetra = tetrahedron_operation(P,P,P,P)
print("Calculated Q_tetra", flush=True)

print("Calculating rhodot", flush=True)
P0 = -(rho*epsilon).bracket(adot[0]).bracket(a2)
P1 = -(rho*epsilon).bracket(a1).bracket(adot[1])
Q_remainder = Q_tetra - P0 - P1
P_withoutrho = -epsilon.bracket(a1).bracket(a2)
rhodot = Q_remainder[0,1] // P_withoutrho[0,1]
print("Calculated rhodot", flush=True)

print("Have nice expression for Q_tetra:", Q_tetra == rhodot*P_withoutrho + P0 + P1, flush=True)

print("Calculating X_rho_formulas", flush=True)
# TODO: Parallelize
X_rho_formulas = [X_formula.bracket(rho*epsilon) for X_formula in X_formulas_independent]
print("Calculated X_rho_formulas", flush=True)

X_rho_basis = set(rhodot.monomials())
for X_rho_formula in X_rho_formulas:
    X_rho_basis |= set(X_rho_formula[0,1,2,3].monomials())
X_rho_basis = list(X_rho_basis)
print("Number of monomials in X_rho_basis:", len(X_rho_basis), flush=True)

print("Calculating X_rho_evaluation_matrix", flush=True)
X_rho_monomial_index = {m : k for k, m in enumerate(X_rho_basis)}
X_rho_evaluation_matrix = matrix(QQ, len(X_rho_basis), len(X_graphs_independent), sparse=True)
for j in range(len(X_graphs_independent)):
    f = X_rho_formulas[j][0,1,2,3]
    v = vector(QQ, len(X_rho_basis), sparse=True)
    for coeff, monomial in zip(f.coefficients(), f.monomials()):
        monomial_index = X_rho_monomial_index[monomial]
        v[monomial_index] = coeff
    X_rho_evaluation_matrix.set_column(j, v)
print("Calculated X_rho_evaluation_matrix", flush=True)

print("Calculating rhodot_vector", flush=True)
rhodot_vector = vector(QQ, len(X_rho_basis), sparse=True)
for coeff, monomial in zip(rhodot.coefficients(), rhodot.monomials()):
    monomial_index = X_rho_monomial_index[monomial]
    rhodot_vector[monomial_index] = coeff
print("Calculated rhodot_vector", flush=True)

# Poisson-triviality $Q_{\text{tetra}}(P) = [\![P, X]\!]$ of the tetrahedral flow

# There is a linear combination of graphs $X$ such that $[\![P, X]\!]$ evaluates to $Q_{\text{tetra}}(P)$:

print("Calculating X_solution_vector", flush=True)
big_matrix = X_a_evaluation_matrix[0].stack(X_a_evaluation_matrix[1]).stack(X_rho_evaluation_matrix)
big_vector = vector(list(-adot_vector[0]) + list(-adot_vector[1]) + list(rhodot_vector))
X_solution_vector = big_matrix.solve_right(big_vector)
print("Calculated X_solution_vector", flush=True)
print("X_solution_vector =", X_solution_vector, flush=True)

X_solution = sum(c*f for c, f in zip(X_solution_vector, X_formulas_independent))

print("P.bracket(X_solution) == Q_tetra:", P.bracket(X_solution) == Q_tetra, flush=True)

# Parameters in the solution $X$ to $[\![P,X]\!] = Q_{\text{tetra}}(P)$

print("Number of parameters in the solution:", big_matrix.right_nullity(), flush=True)

print("Basis of kernel:", big_matrix.right_kernel().basis(), flush=True)
\end{verbatim}
\normalsize
The output:
\tiny
\begin{verbatim}
Number of graphs in X: 160
Calculating X_formulas
Calculated X_formulas
Number of monomials in components of X: [71370, 71370, 71370, 71370]
Calculating X_evaluation_matrix
Calculated X_evaluation_matrix
Number of linearly independent graphs in X: 72
Maximal subset of linearly independent graphs in X: [0, 1, 2, 3, 5, 6, 7, 8, 9, 10, 11, 12, 13, 14, 15, 16, 17, 18, 19, 20, 21, 22, 23,
24, 25, 26, 27, 28, 29, 31, 32, 33, 36, 37, 41, 62, 63, 64, 66, 67, 74, 75, 78, 80, 86, 87, 91, 113, 114, 115, 116, 117, 118, 119, 120,
122, 123, 124, 125, 126, 130, 131, 133, 134, 135, 142, 143, 145, 150, 152, 153, 154]
Calculating adot
Calculated adot
Calculating X_a_formulas
Calculated X_a_formulas
Number of monomials in X_a_basis: 95424 95424
Calculating X_a_evaluation_matrix
Calculated X_a_evaluation_matrix
Calculating adot_vector
Calculated adot_vector
Calculating Q_tetra
Calculated Q_tetra
Calculating rhodot
Calculated rhodot
Have nice expression for Q_tetra: True
Calculating X_rho_formulas
Calculated X_rho_formulas
Number of monomials in X_rho_basis: 311508
Calculating X_rho_evaluation_matrix
Calculated X_rho_evaluation_matrix
Calculating rhodot_vector
Calculated rhodot_vector
Calculating X_solution_vector
Traceback (most recent call last):
  File "/home2/s3058069/sol_3D_9_plus_vanishing_to_4D.sage.py", line 283, in <module>
    X_solution_vector = big_matrix.solve_right(big_vector)
  File "sage/matrix/matrix2.pyx", line 939, in sage.matrix.matrix2.Matrix.solve_right (build/cythonized/sage/matrix/matrix2.c:16151)
  File "sage/matrix/matrix2.pyx", line 1062, in sage.matrix.matrix2.Matrix._solve_right_general (build/cythonized/sage/matrix/matrix2.c:17699)
ValueError: matrix equation has no solutions

###############################################################################
Hábrók Cluster

\end{verbatim}
\normalsize 

\subsection{Solution n.2 over sunflower micro-graphs in 3D does not extend to a 4D solution}\label{A2}

This result is in Ideas \ref{idea1},\ref{idea2}, Propositions \ref{fail1},\ref{fail2}, and lines 3,4,5 and 7,8,9 in Table~\ref{updowntable} on p.~\pageref{updowntable}. This is because in this script, we include the 4D-descendants of the vanishing 3D sunflower micro-graphs directly.

The script:
\tiny
\begin{verbatim}
NPROCS = 32

# Take the #10 over 3D, plus the #13 vanishing (out of full #48 sunflower), extend to 4D: is there a solution?

import warnings
warnings.filterwarnings("ignore", category=DeprecationWarning)

# Graphs for the vector field $X$

X_graph_encodings = [(0,1,4,7,1,3,5,8,1,2,6,9),(0,1,4,7,1,6,5,8,4,2,6,9),(0,1,4,7,1,9,5,8,4,2,6,9),(0,1,4,7,1,6,5,8,7,2,6,9),
(0,1,4,7,1,9,5,8,7,2,6,9),(0,1,4,7,4,3,5,8,4,2,6,9),(0,1,4,7,7,3,5,8,4,2,6,9),(0,1,4,7,7,3,5,8,7,2,6,9),
(0,1,4,7,4,3,5,8,7,2,6,9),(0,2,4,7,1,3,5,8,1,2,6,9),(0,2,4,7,1,3,5,8,1,5,6,9),(0,2,4,7,1,3,5,8,1,8,6,9),
(0,2,4,7,1,6,5,8,1,5,6,9),(0,2,4,7,1,9,5,8,1,5,6,9),(0,2,4,7,1,9,5,8,1,8,6,9),(0,2,4,7,1,6,5,8,1,8,6,9),
(0,2,4,7,4,3,5,8,1,5,6,9),(0,2,4,7,7,3,5,8,1,5,6,9),(0,2,4,7,7,3,5,8,1,8,6,9),(0,2,4,7,4,3,5,8,1,8,6,9),
(0,5,4,7,1,3,5,8,1,2,6,9),(0,8,4,7,1,3,5,8,1,2,6,9)]
for i in [4,7]:
   for j in [6,9]:
       for k in [4,7]:
           for l in [5,8]:
               X_graph_encodings.append((0,1,4,7,i,j,5,8,k,l,6,9))
for i in [4,7]:
    for j in [6,9]:
        for k in [4,7]:
            X_graph_encodings.append((0,1,4,7,i,j,5,8,k,2,6,9))

vanishing_expanded=[(0,2,4,7,1,3,5,8,4,2,6,9),(0,2,4,7,1,3,5,8,7,2,6,9),(0,5,4,7,4,3,5,8,1,2,6,9),(0,8,4,7,4,3,5,8,1,2,6,9),
(0,8,4,7,7,3,5,8,1,2,6,9),(0,5,4,7,7,3,5,8,1,2,6,9),(0,5,4,7,1,3,5,8,4,2,6,9),(0,8,4,7,1,3,5,8,4,2,6,9),
(0,8,4,7,1,3,5,8,7,2,6,9),(0,5,4,7,1,3,5,8,7,2,6,9),(0,2,4,7,1,3,5,8,4,5,6,9),(0,2,4,7,1,3,5,8,7,5,6,9),
(0,2,4,7,1,3,5,8,7,8,6,9),(0,2,4,7,1,3,5,8,4,8,6,9)]

for i in [4,7]:
    for j in [5,8]:
        for k in [6,9]:
            vanishing_expanded.append((0,1,4,7,1,k,5,8,i,j,6,9))
            vanishing_expanded.append((0,2,4,7,1,k,5,8,i,j,6,9))
            vanishing_expanded.append((0,2,4,7,i,k,5,8,1,j,6,9))
            vanishing_expanded.append((0,j,4,7,i,k,5,8,1,2,6,9))
            vanishing_expanded.append((0,j,4,7,1,k,5,8,i,2,6,9))
            vanishing_expanded.append((0,1,4,7,i,k,5,8,1,j,6,9))
            
for i in [4,7]:
    for j1 in [5,8]:
        for j2 in [5,8]:
            vanishing_expanded.append((0,j1,4,7,1,3,5,8,i,j2,6,9))
            
for i1 in [4,7]:
    for i2 in [4,7]:
        for j1 in [5,8]:
            for j2 in [5,8]:
                vanishing_expanded.append((0,j1,4,7,i1,3,5,8,i2,j2,6,9))
                
for i1 in [4,7]:
    for i2 in [4,7]:
        for j1 in [5,8]:
            for j2 in [5,8]:
                for k in [6,9]:
                    vanishing_expanded.append((0,j1,4,7,i1,k,5,8,i2,j2,6,9))

X_graph_encodings += vanishing_expanded

# Convert the encodings to graphs:

from gcaops.graph.formality_graph import FormalityGraph
def encoding_to_graph(encoding):
    targets = [encoding[0:4], encoding[4:8], encoding[8:12]]
    edges = sum([[(k+1,v) for v in t] for (k,t) in enumerate(targets)], [])
    return FormalityGraph(1, 9, edges)

X_graphs = [encoding_to_graph(e) for e in X_graph_encodings]
print("Number of graphs in X:", len(X_graphs), flush=True)

# Formulas of graphs for $X$

# We want to evaluate each graph to a formula, inserting three copies of $\varrho\varepsilon^{ijk}$ 
# and three copies of the Casimir $a$ into the aerial vertices.

# First define the even coordinates $x,y,z$ (base), $\varrho, a$ (fibre) and odd coordinates $\xi_0, \xi_1, \xi_2$:

from gcaops.algebra.differential_polynomial_ring import DifferentialPolynomialRing
D4 = DifferentialPolynomialRing(QQ, ('rho','a1','a2'), ('x','y','z','w'), max_differential_orders=[3+1,1+3+1,1+3+1])
rho, a1, a2 = D4.fibre_variables()
x,y,z,w = D4.base_variables()
even_coords = [x,y,z,w]

from gcaops.algebra.superfunction_algebra import SuperfunctionAlgebra

S4.<xi0,xi1,xi2,xi3> = SuperfunctionAlgebra(D4, D4.base_variables())
xi = S4.gens()
odd_coords = xi
epsilon = xi[0]*xi[1]*xi[2]*xi[3] # Levi-Civita tensor

# Now evaluate each graph to a formula:

import itertools
from multiprocessing import Pool

def evaluate_graph(g):
    E = x*xi[0] + y*xi[1] + z*xi[2] + w*xi[3] # Euler vector field, to insert into ground vertex. Incoming derivative d/dx^i will result in xi[i].
    result = S4.zero()
    for index_choice in itertools.product(itertools.permutations(range(4)), repeat=3):
        sign = epsilon[index_choice[0]] * epsilon[index_choice[1]] * epsilon[index_choice[2]]
        # NOTE: This assumes the ground vertex is labeled 0, the vertices of out-degree 4 are labeled 1, 2, 3, 
        # and the Casimirs are labeled 4, 5, 6; 7, 8, 9.
        vertex_content = [E, S4(rho), S4(rho), S4(rho), S4(a1), S4(a1), S4(a1), S4(a2), S4(a2), S4(a2)]
        for ((source, target), index) in zip(g.edges(), sum(map(list, index_choice), [])):
            vertex_content[target] = vertex_content[target].derivative(even_coords[index])
        result += sign * prod(vertex_content)
    return result

print("Calculating X_formulas", flush=True)
X_formulas = []
with Pool(processes=NPROCS) as pool:
    X_formulas = list(pool.imap(evaluate_graph, X_graphs))
# NOTE: Necessary fixup of parents after multiprocessing.
X_formulas = [S4(X_formula) for X_formula in X_formulas]
print("Calculated X_formulas", flush=True)

# Each resulting formula is of the form $X_g = \sum_{i=1}^3 X_g^i \xi_i$ where $X_g^i$ is a differential polynomial depending on the graph $g$.

# Relations between (formulas of) graphs in $X$

# For each $i = 1,2,3,4$ we collect the distinct differential monomials in $X_g^i$ ranging over all $g$:

X_monomial_basis = [set([]) for i in range(4)]
for i in range(4):
    for X in X_formulas:
        X_monomial_basis[i] |= set(X[i].monomials())
X_monomial_basis = [list(b) for b in X_monomial_basis]
X_monomial_index = [{m : k for k, m in enumerate(b)} for b in X_monomial_basis]

print("Number of monomials in components of X:", [len(b) for b in X_monomial_basis], flush=True)

X_monomial_count = sum(len(b) for b in X_monomial_basis)

# Now the graph-to-formula evaluation can be expressed as a matrix:

print("Calculating X_evaluation_matrix", flush=True)
X_evaluation_matrix = matrix(QQ, X_monomial_count, len(X_graphs), sparse=True)
for i in range(len(X_graphs)):
    v = vector(QQ, X_monomial_count, sparse=True)
    index_shift = 0
    for j in range(4):
        f = X_formulas[i][j]
        for coeff, monomial in zip(f.coefficients(), f.monomials()):
            monomial_index = X_monomial_index[j][monomial]
            v[index_shift + monomial_index] = coeff
        index_shift += len(X_monomial_basis[j])
    X_evaluation_matrix.set_column(i, v)
print("Calculated X_evaluation_matrix", flush=True)

# Here is the number of linearly independent graphs:

print("Number of linearly independent graphs in X:", X_evaluation_matrix.rank(), flush=True)

# We collect the linearly independent graphs and their formulas in new shorter lists:

pivots = X_evaluation_matrix.pivots()
print("Maximal subset of linearly independent graphs in X:", list(pivots), flush=True)
X_graphs_independent = [X_graphs[k] for k in pivots]
X_formulas_independent = [X_formulas[k] for k in pivots]

# The tetrahedron operation $\operatorname{Op}(\gamma_3)$

from gcaops.graph.undirected_graph_complex import UndirectedGraphComplex

GC = UndirectedGraphComplex(QQ, implementation='vector', sparse=True)
tetrahedron = GC.cohomology_basis(4, 6)[0]

from gcaops.graph.directed_graph_complex import DirectedGraphComplex

dGC = DirectedGraphComplex(QQ, implementation='vector')
tetrahedron_oriented = dGC(tetrahedron)
tetrahedron_oriented_filtered = tetrahedron_oriented.filter(max_out_degree=2)
tetrahedron_operation = S4.graph_operation(tetrahedron_oriented_filtered)

# The flow $\dot{a}_i = 4\cdot\operatorname{Op}(\gamma_3)(P,P,P,a_i)$

# Let $P$ be the 4D Nambu–Poisson bracket with pre-factor $\varrho$ and Casimir functions $a_1, a_2$:

P = -(rho*epsilon).bracket(a1).bracket(a2)

# We calculate $\dot{a}_i = 4\cdot\operatorname{Op}(\gamma_3)(P,P,P,a_i)$:

def casimir_flow(f):
    return 4*tetrahedron_operation(P,P,P,f)

a = [S4(a1), S4(a2)]

print("Calculating adot", flush=True)
adot = [casimir_flow(a_i) for a_i in a]
# adot = []
# with Pool(processes=2) as pool:
#    adot = list(pool.imap(casimir_flow, a))
# NOTE: Necessary fixup of parents after multiprocessing.
adot = [S4(D4(adot_i[()])) for adot_i in adot]
print("Calculated adot", flush=True)

# We calculate $X_g(a_i) = [\![X_g, a_i]\!]$ for $i=1,2$:

# TODO: Parallelize
print("Calculating X_a_formulas", flush=True)
X_a_formulas = [[X_formula.bracket(f) for X_formula in X_formulas_independent] for f in a]
print("Calculated X_a_formulas", flush=True)

# We express $g \mapsto X_g(a_i)$ as a sparse matrix for $i=1,2$:

X_a_basis = [set(f[()].monomials()) for f in adot]
for k in range(len(a)):
    for X_a_formula in X_a_formulas[k]:
        X_a_basis[k] |= set(X_a_formula[()].monomials())
X_a_basis = [list(B) for B in X_a_basis]
print("Number of monomials in X_a_basis:", len(X_a_basis[0]), len(X_a_basis[1]), flush=True)

print("Calculating X_a_evaluation_matrix", flush=True)
X_a_monomial_index = [{m : k for k, m in enumerate(B)} for B in X_a_basis]
X_a_evaluation_matrix = [matrix(QQ, len(B), len(X_graphs_independent), sparse=True) for B in X_a_basis]
for i in range(len(a)):
    for j in range(len(X_graphs_independent)):
        v = vector(QQ, len(X_a_basis[i]), sparse=True)
        f = X_a_formulas[i][j][()]
        for coeff, monomial in zip(f.coefficients(), f.monomials()):
            monomial_index = X_a_monomial_index[i][monomial]
            v[monomial_index] = coeff
        X_a_evaluation_matrix[i].set_column(j, v)
print("Calculated X_a_evaluation_matrix", flush=True)

# We express $\dot{a}_i$ as a vector for $i=1,2$:

print("Calculating adot_vector", flush=True)
adot_vector = [vector(QQ, len(B)) for B in X_a_basis]
for i in range(len(a)):
    f = adot[i][()]
    for coeff, monomial in zip(f.coefficients(), f.monomials()):
        monomial_index = X_a_monomial_index[i][monomial]
        adot_vector[i][monomial_index] = coeff
print("Calculated adot_vector", flush=True)

# The flow $\dot{\varrho}$

print("Calculating Q_tetra", flush=True)
Q_tetra = tetrahedron_operation(P,P,P,P)
print("Calculated Q_tetra", flush=True)

print("Calculating rhodot", flush=True)
P0 = -(rho*epsilon).bracket(adot[0]).bracket(a2)
P1 = -(rho*epsilon).bracket(a1).bracket(adot[1])
Q_remainder = Q_tetra - P0 - P1
P_withoutrho = -epsilon.bracket(a1).bracket(a2)
rhodot = Q_remainder[0,1] // P_withoutrho[0,1]
print("Calculated rhodot", flush=True)

print("Have nice expression for Q_tetra:", Q_tetra == rhodot*P_withoutrho + P0 + P1, flush=True)

print("Calculating X_rho_formulas", flush=True)
# TODO: Parallelize
X_rho_formulas = [X_formula.bracket(rho*epsilon) for X_formula in X_formulas_independent]
print("Calculated X_rho_formulas", flush=True)

X_rho_basis = set(rhodot.monomials())
for X_rho_formula in X_rho_formulas:
    X_rho_basis |= set(X_rho_formula[0,1,2,3].monomials())
X_rho_basis = list(X_rho_basis)
print("Number of monomials in X_rho_basis:", len(X_rho_basis), flush=True)

print("Calculating X_rho_evaluation_matrix", flush=True)
X_rho_monomial_index = {m : k for k, m in enumerate(X_rho_basis)}
X_rho_evaluation_matrix = matrix(QQ, len(X_rho_basis), len(X_graphs_independent), sparse=True)
for j in range(len(X_graphs_independent)):
    f = X_rho_formulas[j][0,1,2,3]
    v = vector(QQ, len(X_rho_basis), sparse=True)
    for coeff, monomial in zip(f.coefficients(), f.monomials()):
        monomial_index = X_rho_monomial_index[monomial]
        v[monomial_index] = coeff
    X_rho_evaluation_matrix.set_column(j, v)
print("Calculated X_rho_evaluation_matrix", flush=True)

print("Calculating rhodot_vector", flush=True)
rhodot_vector = vector(QQ, len(X_rho_basis), sparse=True)
for coeff, monomial in zip(rhodot.coefficients(), rhodot.monomials()):
    monomial_index = X_rho_monomial_index[monomial]
    rhodot_vector[monomial_index] = coeff
print("Calculated rhodot_vector", flush=True)

# Poisson-triviality $Q_{\text{tetra}}(P) = [\![P, X]\!]$ of the tetrahedral flow

# There is a linear combination of graphs $X$ such that $[\![P, X]\!]$ evaluates to $Q_{\text{tetra}}(P)$:

print("Calculating X_solution_vector", flush=True)
big_matrix = X_a_evaluation_matrix[0].stack(X_a_evaluation_matrix[1]).stack(X_rho_evaluation_matrix)
big_vector = vector(list(-adot_vector[0]) + list(-adot_vector[1]) + list(rhodot_vector))
X_solution_vector = big_matrix.solve_right(big_vector)
print("Calculated X_solution_vector", flush=True)
print("X_solution_vector =", X_solution_vector, flush=True)

X_solution = sum(c*f for c, f in zip(X_solution_vector, X_formulas_independent))

print("P.bracket(X_solution) == Q_tetra:", P.bracket(X_solution) == Q_tetra, flush=True)

# Parameters in the solution $X$ to $[\![P,X]\!] = Q_{\text{tetra}}(P)$

print("Number of parameters in the solution:", big_matrix.right_nullity(), flush=True)

print("Basis of kernel:", big_matrix.right_kernel().basis(), flush=True)
\end{verbatim}
\normalsize
The output:
\tiny
\begin{verbatim}
Number of graphs in X: 164
Calculating X_formulas
Calculated X_formulas
Number of monomials in components of X: [72588, 72588, 72588, 72588]
Calculating X_evaluation_matrix
Calculated X_evaluation_matrix
Number of linearly independent graphs in X: 75
Maximal subset of linearly independent graphs in X: [0, 1, 2, 3, 4, 5, 6, 7, 9, 10, 11, 12, 13, 14, 16, 17, 18, 19, 20, 21, 22, 23, 24,
25, 27, 28, 29, 32, 33, 37, 38, 39, 40, 41, 42, 43, 44, 45, 66, 67, 68, 70, 71, 78, 79, 82, 84, 90, 91, 95, 117, 118, 119, 120, 121, 122,
123, 124, 126, 127, 128, 129, 130, 134, 135, 137, 138, 139, 146, 147, 149, 154, 156, 157, 158]
Calculating adot
Calculated adot
Calculating X_a_formulas
Calculated X_a_formulas
Number of monomials in X_a_basis: 103944 103944
Calculating X_a_evaluation_matrix
Calculated X_a_evaluation_matrix
Calculating adot_vector
Calculated adot_vector
Calculating Q_tetra
Calculated Q_tetra
Calculating rhodot
Calculated rhodot
Have nice expression for Q_tetra: True
Calculating X_rho_formulas
Calculated X_rho_formulas
Number of monomials in X_rho_basis: 318420
Calculating X_rho_evaluation_matrix
Calculated X_rho_evaluation_matrix
Calculating rhodot_vector
Calculated rhodot_vector
Calculating X_solution_vector
Traceback (most recent call last):
  File "/home2/s3058069/sol_3D_10_plus_vanishing_to_4D.sage.py", line 286, in <module>
    X_solution_vector = big_matrix.solve_right(big_vector)
  File "sage/matrix/matrix2.pyx", line 939, in sage.matrix.matrix2.Matrix.solve_right (build/cythonized/sage/matrix/matrix2.c:16151)
  File "sage/matrix/matrix2.pyx", line 1062, in sage.matrix.matrix2.Matrix._solve_right_general (build/cythonized/sage/matrix/matrix2.c:17699)
ValueError: matrix equation has no solutions

###############################################################################
Hábrók Cluster

\end{verbatim}
\normalsize 

\subsection{Solution n.3 over sunflower micro-graphs in 3D does not extend to a 4D solution}\label{A3}

This result is in Ideas \ref{idea1},\ref{idea2}, Propositions \ref{fail1},\ref{fail2}, and lines 3,4,5 and 7,8,9 in Table~\ref{updowntable} on p.~\pageref{updowntable}. This is because in this script, we include the 4D-descendants of the vanishing 3D sunflower micro-graphs directly.

The script:
\tiny
\begin{verbatim}
NPROCS = 32

# Take the #12 over 3D, from 4D sol (#56) projected down as linear combination, plus the #13 vanishing (out of full #48 sunflower), extend to 4D:
# is there a solution?

import warnings
warnings.filterwarnings("ignore", category=DeprecationWarning)

# Graphs for the vector field $X$

X_graph_encodings = [(0,1,4,7,1,3,5,8,1,2,6,9),(0,1,4,7,1,6,5,8,4,2,6,9),(0,1,4,7,1,9,5,8,4,2,6,9),(0,1,4,7,1,6,5,8,7,2,6,9),
(0,1,4,7,1,9,5,8,7,2,6,9),(0,1,4,7,4,3,5,8,4,2,6,9),(0,1,4,7,7,3,5,8,4,2,6,9),(0,1,4,7,7,3,5,8,7,2,6,9),
(0,1,4,7,4,3,5,8,7,2,6,9),(0,2,4,7,1,3,5,8,1,2,6,9),(0,2,4,7,1,3,5,8,1,5,6,9),(0,2,4,7,1,3,5,8,1,8,6,9),
(0,2,4,7,1,6,5,8,1,5,6,9),(0,2,4,7,1,9,5,8,1,5,6,9),(0,2,4,7,1,9,5,8,1,8,6,9),(0,2,4,7,1,6,5,8,1,8,6,9),
(0,2,4,7,4,3,5,8,1,5,6,9),(0,2,4,7,7,3,5,8,1,5,6,9),(0,2,4,7,7,3,5,8,1,8,6,9),(0,2,4,7,4,3,5,8,1,8,6,9),
(0,5,4,7,1,3,5,8,1,2,6,9),(0,8,4,7,1,3,5,8,1,2,6,9),(0,1,4,7,1,6,5,8,1,5,6,9),(0,1,4,7,1,9,5,8,1,5,6,9),
(0,1,4,7,1,9,5,8,1,8,6,9),(0,1,4,7,1,6,5,8,1,8,6,9),(0,5,4,7,1,6,5,8,1,5,6,9),(0,5,4,7,1,6,5,8,1,8,6,9),
(0,5,4,7,1,9,5,8,1,5,6,9),(0,5,4,7,1,9,5,8,1,8,6,9),(0,8,4,7,1,6,5,8,1,5,6,9),(0,8,4,7,1,6,5,8,1,8,6,9),
(0,8,4,7,1,9,5,8,1,5,6,9),(0,8,4,7,1,9,5,8,1,8,6,9)]
for i in [4,7]:
   for j in [6,9]:
       for k in [4,7]:
           for l in [5,8]:
               X_graph_encodings.append((0,1,4,7,i,j,5,8,k,l,6,9))
for i in [4,7]:
    for j in [6,9]:
        for k in [4,7]:
            X_graph_encodings.append((0,1,4,7,i,j,5,8,k,2,6,9))

vanishing_expanded=[(0,2,4,7,1,3,5,8,4,2,6,9),(0,2,4,7,1,3,5,8,7,2,6,9),(0,5,4,7,4,3,5,8,1,2,6,9),(0,8,4,7,4,3,5,8,1,2,6,9),
(0,8,4,7,7,3,5,8,1,2,6,9),(0,5,4,7,7,3,5,8,1,2,6,9),(0,5,4,7,1,3,5,8,4,2,6,9),(0,8,4,7,1,3,5,8,4,2,6,9),
(0,8,4,7,1,3,5,8,7,2,6,9),(0,5,4,7,1,3,5,8,7,2,6,9),(0,2,4,7,1,3,5,8,4,5,6,9),(0,2,4,7,1,3,5,8,7,5,6,9),
(0,2,4,7,1,3,5,8,7,8,6,9),(0,2,4,7,1,3,5,8,4,8,6,9)]

for i in [4,7]:
    for j in [5,8]:
        for k in [6,9]:
            vanishing_expanded.append((0,1,4,7,1,k,5,8,i,j,6,9))
            vanishing_expanded.append((0,2,4,7,1,k,5,8,i,j,6,9))
            vanishing_expanded.append((0,2,4,7,i,k,5,8,1,j,6,9))
            vanishing_expanded.append((0,j,4,7,i,k,5,8,1,2,6,9))
            vanishing_expanded.append((0,j,4,7,1,k,5,8,i,2,6,9))
            vanishing_expanded.append((0,1,4,7,i,k,5,8,1,j,6,9))
            
for i in [4,7]:
    for j1 in [5,8]:
        for j2 in [5,8]:
            vanishing_expanded.append((0,j1,4,7,1,3,5,8,i,j2,6,9))
            
for i1 in [4,7]:
    for i2 in [4,7]:
        for j1 in [5,8]:
            for j2 in [5,8]:
                vanishing_expanded.append((0,j1,4,7,i1,3,5,8,i2,j2,6,9))
                
for i1 in [4,7]:
    for i2 in [4,7]:
        for j1 in [5,8]:
            for j2 in [5,8]:
                for k in [6,9]:
                    vanishing_expanded.append((0,j1,4,7,i1,k,5,8,i2,j2,6,9))

X_graph_encodings += vanishing_expanded


# Convert the encodings to graphs:

from gcaops.graph.formality_graph import FormalityGraph
def encoding_to_graph(encoding):
    targets = [encoding[0:4], encoding[4:8], encoding[8:12]]
    edges = sum([[(k+1,v) for v in t] for (k,t) in enumerate(targets)], [])
    return FormalityGraph(1, 9, edges)

X_graphs = [encoding_to_graph(e) for e in X_graph_encodings]
print("Number of graphs in X:", len(X_graphs), flush=True)

# Formulas of graphs for $X$

# We want to evaluate each graph to a formula, inserting three copies of $\varrho\varepsilon^{ijk}$ 
# and three copies of the Casimir $a$ into the aerial vertices.

# First define the even coordinates $x,y,z$ (base), $\varrho, a$ (fibre) and odd coordinates $\xi_0, \xi_1, \xi_2$:

from gcaops.algebra.differential_polynomial_ring import DifferentialPolynomialRing
D4 = DifferentialPolynomialRing(QQ, ('rho','a1','a2'), ('x','y','z','w'), max_differential_orders=[3+1,1+3+1,1+3+1])
rho, a1, a2 = D4.fibre_variables()
x,y,z,w = D4.base_variables()
even_coords = [x,y,z,w]

from gcaops.algebra.superfunction_algebra import SuperfunctionAlgebra

S4.<xi0,xi1,xi2,xi3> = SuperfunctionAlgebra(D4, D4.base_variables())
xi = S4.gens()
odd_coords = xi
epsilon = xi[0]*xi[1]*xi[2]*xi[3] # Levi-Civita tensor

# Now evaluate each graph to a formula:

import itertools
from multiprocessing import Pool

def evaluate_graph(g):
    E = x*xi[0] + y*xi[1] + z*xi[2] + w*xi[3] # Euler vector field, to insert into ground vertex. Incoming derivative d/dx^i will result in xi[i].
    result = S4.zero()
    for index_choice in itertools.product(itertools.permutations(range(4)), repeat=3):
        sign = epsilon[index_choice[0]] * epsilon[index_choice[1]] * epsilon[index_choice[2]]
        # NOTE: This assumes the ground vertex is labeled 0, the vertices of out-degree 4 are labeled 1, 2, 3, 
        # and the Casimirs are labeled 4, 5, 6; 7, 8, 9.
        vertex_content = [E, S4(rho), S4(rho), S4(rho), S4(a1), S4(a1), S4(a1), S4(a2), S4(a2), S4(a2)]
        for ((source, target), index) in zip(g.edges(), sum(map(list, index_choice), [])):
            vertex_content[target] = vertex_content[target].derivative(even_coords[index])
        result += sign * prod(vertex_content)
    return result

print("Calculating X_formulas", flush=True)
X_formulas = []
with Pool(processes=NPROCS) as pool:
    X_formulas = list(pool.imap(evaluate_graph, X_graphs))
# NOTE: Necessary fixup of parents after multiprocessing.
X_formulas = [S4(X_formula) for X_formula in X_formulas]
print("Calculated X_formulas", flush=True)

# Each resulting formula is of the form $X_g = \sum_{i=1}^3 X_g^i \xi_i$ where $X_g^i$ is a differential polynomial depending on the graph $g$.

# Relations between (formulas of) graphs in $X$

# For each $i = 1,2,3,4$ we collect the distinct differential monomials in $X_g^i$ ranging over all $g$:

X_monomial_basis = [set([]) for i in range(4)]
for i in range(4):
    for X in X_formulas:
        X_monomial_basis[i] |= set(X[i].monomials())
X_monomial_basis = [list(b) for b in X_monomial_basis]
X_monomial_index = [{m : k for k, m in enumerate(b)} for b in X_monomial_basis]

print("Number of monomials in components of X:", [len(b) for b in X_monomial_basis], flush=True)

X_monomial_count = sum(len(b) for b in X_monomial_basis)

# Now the graph-to-formula evaluation can be expressed as a matrix:

print("Calculating X_evaluation_matrix", flush=True)
X_evaluation_matrix = matrix(QQ, X_monomial_count, len(X_graphs), sparse=True)
for i in range(len(X_graphs)):
    v = vector(QQ, X_monomial_count, sparse=True)
    index_shift = 0
    for j in range(4):
        f = X_formulas[i][j]
        for coeff, monomial in zip(f.coefficients(), f.monomials()):
            monomial_index = X_monomial_index[j][monomial]
            v[index_shift + monomial_index] = coeff
        index_shift += len(X_monomial_basis[j])
    X_evaluation_matrix.set_column(i, v)
print("Calculated X_evaluation_matrix", flush=True)

# Here is the number of linearly independent graphs:

print("Number of linearly independent graphs in X:", X_evaluation_matrix.rank(), flush=True)

# We collect the linearly independent graphs and their formulas in new shorter lists:

pivots = X_evaluation_matrix.pivots()
print("Maximal subset of linearly independent graphs in X:", list(pivots), flush=True)
X_graphs_independent = [X_graphs[k] for k in pivots]
X_formulas_independent = [X_formulas[k] for k in pivots]

# The tetrahedron operation $\operatorname{Op}(\gamma_3)$

from gcaops.graph.undirected_graph_complex import UndirectedGraphComplex

GC = UndirectedGraphComplex(QQ, implementation='vector', sparse=True)
tetrahedron = GC.cohomology_basis(4, 6)[0]

from gcaops.graph.directed_graph_complex import DirectedGraphComplex

dGC = DirectedGraphComplex(QQ, implementation='vector')
tetrahedron_oriented = dGC(tetrahedron)
tetrahedron_oriented_filtered = tetrahedron_oriented.filter(max_out_degree=2)
tetrahedron_operation = S4.graph_operation(tetrahedron_oriented_filtered)

# The flow $\dot{a}_i = 4\cdot\operatorname{Op}(\gamma_3)(P,P,P,a_i)$

# Let $P$ be the 4D Nambu–Poisson bracket with pre-factor $\varrho$ and Casimir functions $a_1, a_2$:

P = -(rho*epsilon).bracket(a1).bracket(a2)

# We calculate $\dot{a}_i = 4\cdot\operatorname{Op}(\gamma_3)(P,P,P,a_i)$:

def casimir_flow(f):
    return 4*tetrahedron_operation(P,P,P,f)

a = [S4(a1), S4(a2)]

print("Calculating adot", flush=True)
adot = [casimir_flow(a_i) for a_i in a]
# adot = []
# with Pool(processes=2) as pool:
#    adot = list(pool.imap(casimir_flow, a))
# NOTE: Necessary fixup of parents after multiprocessing.
adot = [S4(D4(adot_i[()])) for adot_i in adot]
print("Calculated adot", flush=True)

# We calculate $X_g(a_i) = [\![X_g, a_i]\!]$ for $i=1,2$:

# TODO: Parallelize
print("Calculating X_a_formulas", flush=True)
X_a_formulas = [[X_formula.bracket(f) for X_formula in X_formulas_independent] for f in a]
print("Calculated X_a_formulas", flush=True)

# We express $g \mapsto X_g(a_i)$ as a sparse matrix for $i=1,2$:

X_a_basis = [set(f[()].monomials()) for f in adot]
for k in range(len(a)):
    for X_a_formula in X_a_formulas[k]:
        X_a_basis[k] |= set(X_a_formula[()].monomials())
X_a_basis = [list(B) for B in X_a_basis]
print("Number of monomials in X_a_basis:", len(X_a_basis[0]), len(X_a_basis[1]), flush=True)

print("Calculating X_a_evaluation_matrix", flush=True)
X_a_monomial_index = [{m : k for k, m in enumerate(B)} for B in X_a_basis]
X_a_evaluation_matrix = [matrix(QQ, len(B), len(X_graphs_independent), sparse=True) for B in X_a_basis]
for i in range(len(a)):
    for j in range(len(X_graphs_independent)):
        v = vector(QQ, len(X_a_basis[i]), sparse=True)
        f = X_a_formulas[i][j][()]
        for coeff, monomial in zip(f.coefficients(), f.monomials()):
            monomial_index = X_a_monomial_index[i][monomial]
            v[monomial_index] = coeff
        X_a_evaluation_matrix[i].set_column(j, v)
print("Calculated X_a_evaluation_matrix", flush=True)

# We express $\dot{a}_i$ as a vector for $i=1,2$:

print("Calculating adot_vector", flush=True)
adot_vector = [vector(QQ, len(B)) for B in X_a_basis]
for i in range(len(a)):
    f = adot[i][()]
    for coeff, monomial in zip(f.coefficients(), f.monomials()):
        monomial_index = X_a_monomial_index[i][monomial]
        adot_vector[i][monomial_index] = coeff
print("Calculated adot_vector", flush=True)

# The flow $\dot{\varrho}$

print("Calculating Q_tetra", flush=True)
Q_tetra = tetrahedron_operation(P,P,P,P)
print("Calculated Q_tetra", flush=True)

print("Calculating rhodot", flush=True)
P0 = -(rho*epsilon).bracket(adot[0]).bracket(a2)
P1 = -(rho*epsilon).bracket(a1).bracket(adot[1])
Q_remainder = Q_tetra - P0 - P1
P_withoutrho = -epsilon.bracket(a1).bracket(a2)
rhodot = Q_remainder[0,1] // P_withoutrho[0,1]
print("Calculated rhodot", flush=True)

print("Have nice expression for Q_tetra:", Q_tetra == rhodot*P_withoutrho + P0 + P1, flush=True)

print("Calculating X_rho_formulas", flush=True)
# TODO: Parallelize
X_rho_formulas = [X_formula.bracket(rho*epsilon) for X_formula in X_formulas_independent]
print("Calculated X_rho_formulas", flush=True)

X_rho_basis = set(rhodot.monomials())
for X_rho_formula in X_rho_formulas:
    X_rho_basis |= set(X_rho_formula[0,1,2,3].monomials())
X_rho_basis = list(X_rho_basis)
print("Number of monomials in X_rho_basis:", len(X_rho_basis), flush=True)

print("Calculating X_rho_evaluation_matrix", flush=True)
X_rho_monomial_index = {m : k for k, m in enumerate(X_rho_basis)}
X_rho_evaluation_matrix = matrix(QQ, len(X_rho_basis), len(X_graphs_independent), sparse=True)
for j in range(len(X_graphs_independent)):
    f = X_rho_formulas[j][0,1,2,3]
    v = vector(QQ, len(X_rho_basis), sparse=True)
    for coeff, monomial in zip(f.coefficients(), f.monomials()):
        monomial_index = X_rho_monomial_index[monomial]
        v[monomial_index] = coeff
    X_rho_evaluation_matrix.set_column(j, v)
print("Calculated X_rho_evaluation_matrix", flush=True)

print("Calculating rhodot_vector", flush=True)
rhodot_vector = vector(QQ, len(X_rho_basis), sparse=True)
for coeff, monomial in zip(rhodot.coefficients(), rhodot.monomials()):
    monomial_index = X_rho_monomial_index[monomial]
    rhodot_vector[monomial_index] = coeff
print("Calculated rhodot_vector", flush=True)

# Poisson-triviality $Q_{\text{tetra}}(P) = [\![P, X]\!]$ of the tetrahedral flow

# There is a linear combination of graphs $X$ such that $[\![P, X]\!]$ evaluates to $Q_{\text{tetra}}(P)$:

print("Calculating X_solution_vector", flush=True)
big_matrix = X_a_evaluation_matrix[0].stack(X_a_evaluation_matrix[1]).stack(X_rho_evaluation_matrix)
big_vector = vector(list(-adot_vector[0]) + list(-adot_vector[1]) + list(rhodot_vector))
X_solution_vector = big_matrix.solve_right(big_vector)
print("Calculated X_solution_vector", flush=True)
print("X_solution_vector =", X_solution_vector, flush=True)

X_solution = sum(c*f for c, f in zip(X_solution_vector, X_formulas_independent))

print("P.bracket(X_solution) == Q_tetra:", P.bracket(X_solution) == Q_tetra, flush=True)

# Parameters in the solution $X$ to $[\![P,X]\!] = Q_{\text{tetra}}(P)$

print("Number of parameters in the solution:", big_matrix.right_nullity(), flush=True)

print("Basis of kernel:", big_matrix.right_kernel().basis(), flush=True)
\end{verbatim}
\normalsize
The output:
\tiny
\begin{verbatim}
Number of graphs in X: 176
Calculating X_formulas
Calculated X_formulas
Number of monomials in components of X: [83946, 83946, 83946, 83946]
Calculating X_evaluation_matrix
Calculated X_evaluation_matrix
Number of linearly independent graphs in X: 86
Maximal subset of linearly independent graphs in X: [0, 1, 2, 3, 4, 5, 6, 7, 9, 10, 11, 12, 13, 14, 16, 17, 18, 19, 20, 21, 22, 23, 24, 26,
27, 28, 29, 30, 31, 32, 33, 34, 35, 36, 37, 39, 40, 41, 44, 45, 49, 50, 51, 52, 53, 54, 55, 56, 57, 78, 79, 80, 82, 83, 90, 91, 94, 96, 102,
103, 107, 129, 130, 131, 132, 133, 134, 135, 136, 138, 139, 140, 141, 142, 146, 147, 149, 150, 151, 158, 159, 161, 166, 168, 169, 170]
Calculating adot
Calculated adot
Calculating X_a_formulas
Calculated X_a_formulas
Number of monomials in X_a_basis: 103944 103944
Calculating X_a_evaluation_matrix
Calculated X_a_evaluation_matrix
Calculating adot_vector
Calculated adot_vector
Calculating Q_tetra
Calculated Q_tetra
Calculating rhodot
Calculated rhodot
Have nice expression for Q_tetra: True
Calculating X_rho_formulas
Calculated X_rho_formulas
Number of monomials in X_rho_basis: 334164
Calculating X_rho_evaluation_matrix
Calculated X_rho_evaluation_matrix
Calculating rhodot_vector
Calculated rhodot_vector
Calculating X_solution_vector
Traceback (most recent call last):
  File "/home2/s3058069/sol_3D_12_plus_vanishing_to_4D.sage.py", line 287, in <module>
    X_solution_vector = big_matrix.solve_right(big_vector)
  File "sage/matrix/matrix2.pyx", line 939, in sage.matrix.matrix2.Matrix.solve_right (build/cythonized/sage/matrix/matrix2.c:16151)
  File "sage/matrix/matrix2.pyx", line 1062, in sage.matrix.matrix2.Matrix._solve_right_general (build/cythonized/sage/matrix/matrix2.c:17699)
ValueError: matrix equation has no solutions

###############################################################################
Hábrók Cluster

\end{verbatim}
\normalsize 

\subsection{The projected 4D to 3D solution $\vec{X}^{\gamma_3}_{d=4}(P)$ under $a^2=w$ does not extend to a 4D solution}\label{A4}

This result is in Idea \ref{4dto3d}, Proposition \ref{fail3} and lines 10,11 in Table~\ref{updowntable} on p.~\pageref{updowntable}. 

The script:
\tiny
\begin{verbatim}
NPROCS = 32

# Take the #12 over 3D which are projected down from the skewed 4D solution, extend to 4D: is there a solution?

# Also add the vanishing-expanded

import warnings
warnings.filterwarnings("ignore", category=DeprecationWarning)

X_graph_encodings = [(0,1,4,7,1,3,5,8,1,2,6,9),(0,2,4,7,1,3,5,8,1,2,6,9)]

for i in [5,8]:
    X_graph_encodings.append((0,2,4,7,1,3,5,8,1,i,6,9))
    X_graph_encodings.append((0,i,4,7,1,3,5,8,1,2,6,9))
    
for i in [4,7]:
    for j in [5,8]:
        X_graph_encodings.append((0,2,4,7,i,3,5,8,1,j,6,9))
        
for i in [5,8]:
    for j in [6,9]:
        for k in [5,8]:
            X_graph_encodings.append((0,i,4,7,1,j,5,8,1,k,6,9))
    
for i in [4,7]:
    for j in [6,9]:
        X_graph_encodings.append((0,1,4,7,1,j,5,8,i,2,6,9))
        
for i in [4,7]:
    for j in [4,7]:
        X_graph_encodings.append((0,1,4,7,i,3,5,8,j,2,6,9))
        
for i in [4,7]:
    for j in [4,7]:
        for k in [6,9]:
            X_graph_encodings.append((0,1,4,7,i,k,5,8,j,2,6,9))
            
for i in [6,9]:
    for j in [5,8]:
        X_graph_encodings.append((0,1,4,7,1,i,5,8,1,j,6,9))
        X_graph_encodings.append((0,2,4,7,1,i,5,8,1,j,6,9))
        
for i in [4,7]:
    for j in [6,9]:
        for k in [4,7]:
            for l in [5,8]:
                X_graph_encodings.append((0,1,4,7,i,j,5,8,k,l,6,9))

vanishing_expanded=[(0,2,4,7,1,3,5,8,4,2,6,9),(0,2,4,7,1,3,5,8,7,2,6,9),(0,5,4,7,4,3,5,8,1,2,6,9),(0,8,4,7,4,3,5,8,1,2,6,9),
(0,8,4,7,7,3,5,8,1,2,6,9),(0,5,4,7,7,3,5,8,1,2,6,9),(0,5,4,7,1,3,5,8,4,2,6,9),(0,8,4,7,1,3,5,8,4,2,6,9),
(0,8,4,7,1,3,5,8,7,2,6,9),(0,5,4,7,1,3,5,8,7,2,6,9),(0,2,4,7,1,3,5,8,4,5,6,9),(0,2,4,7,1,3,5,8,7,5,6,9),
(0,2,4,7,1,3,5,8,7,8,6,9),(0,2,4,7,1,3,5,8,4,8,6,9)]

for i in [4,7]:
    for j in [5,8]:
        for k in [6,9]:
            vanishing_expanded.append((0,1,4,7,1,k,5,8,i,j,6,9))
            vanishing_expanded.append((0,2,4,7,1,k,5,8,i,j,6,9))
            vanishing_expanded.append((0,2,4,7,i,k,5,8,1,j,6,9))
            vanishing_expanded.append((0,j,4,7,i,k,5,8,1,2,6,9))
            vanishing_expanded.append((0,j,4,7,1,k,5,8,i,2,6,9))
            vanishing_expanded.append((0,1,4,7,i,k,5,8,1,j,6,9))
            
for i in [4,7]:
    for j1 in [5,8]:
        for j2 in [5,8]:
            vanishing_expanded.append((0,j1,4,7,1,3,5,8,i,j2,6,9))
            
for i1 in [4,7]:
    for i2 in [4,7]:
        for j1 in [5,8]:
            for j2 in [5,8]:
                vanishing_expanded.append((0,j1,4,7,i1,3,5,8,i2,j2,6,9))
                
for i1 in [4,7]:
    for i2 in [4,7]:
        for j1 in [5,8]:
            for j2 in [5,8]:
                for k in [6,9]:
                    vanishing_expanded.append((0,j1,4,7,i1,k,5,8,i2,j2,6,9))

X_graph_encodings += vanishing_expanded

# Convert the encodings to graphs:

from gcaops.graph.formality_graph import FormalityGraph
def encoding_to_graph(encoding):
    targets = [encoding[0:4], encoding[4:8], encoding[8:12]]
    edges = sum([[(k+1,v) for v in t] for (k,t) in enumerate(targets)], [])
    return FormalityGraph(1, 9, edges)

X_graphs = [encoding_to_graph(e) for e in X_graph_encodings]
print("Number of graphs in X:", len(X_graphs), flush=True)

# Formulas of graphs for $X$

# We want to evaluate each graph to a formula, inserting three copies of $\varrho\varepsilon^{ijk}$ 
# and three copies of the Casimir $a$ into the aerial vertices.

# First define the even coordinates $x,y,z$ (base), $\varrho, a$ (fibre) and odd coordinates $\xi_0, \xi_1, \xi_2$:

from gcaops.algebra.differential_polynomial_ring import DifferentialPolynomialRing
D4 = DifferentialPolynomialRing(QQ, ('rho','a1','a2'), ('x','y','z','w'), max_differential_orders=[3+1,1+3+1,1+3+1])
rho, a1, a2 = D4.fibre_variables()
x,y,z,w = D4.base_variables()
even_coords = [x,y,z,w]

from gcaops.algebra.superfunction_algebra import SuperfunctionAlgebra

S4.<xi0,xi1,xi2,xi3> = SuperfunctionAlgebra(D4, D4.base_variables())
xi = S4.gens()
odd_coords = xi
epsilon = xi[0]*xi[1]*xi[2]*xi[3] # Levi-Civita tensor

# Now evaluate each graph to a formula:

import itertools
from multiprocessing import Pool

def evaluate_graph(g):
    E = x*xi[0] + y*xi[1] + z*xi[2] + w*xi[3] # Euler vector field, to insert into ground vertex. Incoming derivative d/dx^i will result in xi[i].
    result = S4.zero()
    for index_choice in itertools.product(itertools.permutations(range(4)), repeat=3):
        sign = epsilon[index_choice[0]] * epsilon[index_choice[1]] * epsilon[index_choice[2]]
        # NOTE: This assumes the ground vertex is labeled 0, the vertices of out-degree 4 are labeled 1, 2, 3, 
        # and the Casimirs are labeled 4, 5, 6; 7, 8, 9.
        vertex_content = [E, S4(rho), S4(rho), S4(rho), S4(a1), S4(a1), S4(a1), S4(a2), S4(a2), S4(a2)]
        for ((source, target), index) in zip(g.edges(), sum(map(list, index_choice), [])):
            vertex_content[target] = vertex_content[target].derivative(even_coords[index])
        result += sign * prod(vertex_content)
    return result

print("Calculating X_formulas", flush=True)
X_formulas = []
with Pool(processes=NPROCS) as pool:
    X_formulas = list(pool.imap(evaluate_graph, X_graphs))
# NOTE: Necessary fixup of parents after multiprocessing.
X_formulas = [S4(X_formula) for X_formula in X_formulas]
print("Calculated X_formulas", flush=True)

# Each resulting formula is of the form $X_g = \sum_{i=1}^3 X_g^i \xi_i$ where $X_g^i$ is a differential polynomial depending on the graph $g$.

# Relations between (formulas of) graphs in $X$

# For each $i = 1,2,3,4$ we collect the distinct differential monomials in $X_g^i$ ranging over all $g$:

X_monomial_basis = [set([]) for i in range(4)]
for i in range(4):
    for X in X_formulas:
        X_monomial_basis[i] |= set(X[i].monomials())
X_monomial_basis = [list(b) for b in X_monomial_basis]
X_monomial_index = [{m : k for k, m in enumerate(b)} for b in X_monomial_basis]

print("Number of monomials in components of X:", [len(b) for b in X_monomial_basis], flush=True)

X_monomial_count = sum(len(b) for b in X_monomial_basis)

# Now the graph-to-formula evaluation can be expressed as a matrix:

print("Calculating X_evaluation_matrix", flush=True)
X_evaluation_matrix = matrix(QQ, X_monomial_count, len(X_graphs), sparse=True)
for i in range(len(X_graphs)):
    v = vector(QQ, X_monomial_count, sparse=True)
    index_shift = 0
    for j in range(4):
        f = X_formulas[i][j]
        for coeff, monomial in zip(f.coefficients(), f.monomials()):
            monomial_index = X_monomial_index[j][monomial]
            v[index_shift + monomial_index] = coeff
        index_shift += len(X_monomial_basis[j])
    X_evaluation_matrix.set_column(i, v)
print("Calculated X_evaluation_matrix", flush=True)

# Here is the number of linearly independent graphs:

print("Number of linearly independent graphs in X:", X_evaluation_matrix.rank(), flush=True)

# We collect the linearly independent graphs and their formulas in new shorter lists:

pivots = X_evaluation_matrix.pivots()
print("Maximal subset of linearly independent graphs in X:", list(pivots), flush=True)
X_graphs_independent = [X_graphs[k] for k in pivots]
X_formulas_independent = [X_formulas[k] for k in pivots]

# The tetrahedron operation $\operatorname{Op}(\gamma_3)$

from gcaops.graph.undirected_graph_complex import UndirectedGraphComplex

GC = UndirectedGraphComplex(QQ, implementation='vector', sparse=True)
tetrahedron = GC.cohomology_basis(4, 6)[0]

from gcaops.graph.directed_graph_complex import DirectedGraphComplex

dGC = DirectedGraphComplex(QQ, implementation='vector')
tetrahedron_oriented = dGC(tetrahedron)
tetrahedron_oriented_filtered = tetrahedron_oriented.filter(max_out_degree=2)
tetrahedron_operation = S4.graph_operation(tetrahedron_oriented_filtered)

# The flow $\dot{a}_i = 4\cdot\operatorname{Op}(\gamma_3)(P,P,P,a_i)$

# Let $P$ be the 4D Nambu–Poisson bracket with pre-factor $\varrho$ and Casimir functions $a_1, a_2$:

P = -(rho*epsilon).bracket(a1).bracket(a2)

# We calculate $\dot{a}_i = 4\cdot\operatorname{Op}(\gamma_3)(P,P,P,a_i)$:

def casimir_flow(f):
    return 4*tetrahedron_operation(P,P,P,f)

a = [S4(a1), S4(a2)]

print("Calculating adot", flush=True)
adot = [casimir_flow(a_i) for a_i in a]
# adot = []
# with Pool(processes=2) as pool:
#    adot = list(pool.imap(casimir_flow, a))
# NOTE: Necessary fixup of parents after multiprocessing.
adot = [S4(D4(adot_i[()])) for adot_i in adot]
print("Calculated adot", flush=True)

# We calculate $X_g(a_i) = [\![X_g, a_i]\!]$ for $i=1,2$:

# TODO: Parallelize
print("Calculating X_a_formulas", flush=True)
X_a_formulas = [[X_formula.bracket(f) for X_formula in X_formulas_independent] for f in a]
print("Calculated X_a_formulas", flush=True)

# We express $g \mapsto X_g(a_i)$ as a sparse matrix for $i=1,2$:

X_a_basis = [set(f[()].monomials()) for f in adot]
for k in range(len(a)):
    for X_a_formula in X_a_formulas[k]:
        X_a_basis[k] |= set(X_a_formula[()].monomials())
X_a_basis = [list(B) for B in X_a_basis]
print("Number of monomials in X_a_basis:", len(X_a_basis[0]), len(X_a_basis[1]), flush=True)

print("Calculating X_a_evaluation_matrix", flush=True)
X_a_monomial_index = [{m : k for k, m in enumerate(B)} for B in X_a_basis]
X_a_evaluation_matrix = [matrix(QQ, len(B), len(X_graphs_independent), sparse=True) for B in X_a_basis]
for i in range(len(a)):
    for j in range(len(X_graphs_independent)):
        v = vector(QQ, len(X_a_basis[i]), sparse=True)
        f = X_a_formulas[i][j][()]
        for coeff, monomial in zip(f.coefficients(), f.monomials()):
            monomial_index = X_a_monomial_index[i][monomial]
            v[monomial_index] = coeff
        X_a_evaluation_matrix[i].set_column(j, v)
print("Calculated X_a_evaluation_matrix", flush=True)

# We express $\dot{a}_i$ as a vector for $i=1,2$:

print("Calculating adot_vector", flush=True)
adot_vector = [vector(QQ, len(B)) for B in X_a_basis]
for i in range(len(a)):
    f = adot[i][()]
    for coeff, monomial in zip(f.coefficients(), f.monomials()):
        monomial_index = X_a_monomial_index[i][monomial]
        adot_vector[i][monomial_index] = coeff
print("Calculated adot_vector", flush=True)

# The flow $\dot{\varrho}$

print("Calculating Q_tetra", flush=True)
Q_tetra = tetrahedron_operation(P,P,P,P)
print("Calculated Q_tetra", flush=True)

print("Calculating rhodot", flush=True)
P0 = -(rho*epsilon).bracket(adot[0]).bracket(a2)
P1 = -(rho*epsilon).bracket(a1).bracket(adot[1])
Q_remainder = Q_tetra - P0 - P1
P_withoutrho = -epsilon.bracket(a1).bracket(a2)
rhodot = Q_remainder[0,1] // P_withoutrho[0,1]
print("Calculated rhodot", flush=True)

print("Have nice expression for Q_tetra:", Q_tetra == rhodot*P_withoutrho + P0 + P1, flush=True)

print("Calculating X_rho_formulas", flush=True)
# TODO: Parallelize
X_rho_formulas = [X_formula.bracket(rho*epsilon) for X_formula in X_formulas_independent]
print("Calculated X_rho_formulas", flush=True)

X_rho_basis = set(rhodot.monomials())
for X_rho_formula in X_rho_formulas:
    X_rho_basis |= set(X_rho_formula[0,1,2,3].monomials())
X_rho_basis = list(X_rho_basis)
print("Number of monomials in X_rho_basis:", len(X_rho_basis), flush=True)

print("Calculating X_rho_evaluation_matrix", flush=True)
X_rho_monomial_index = {m : k for k, m in enumerate(X_rho_basis)}
X_rho_evaluation_matrix = matrix(QQ, len(X_rho_basis), len(X_graphs_independent), sparse=True)
for j in range(len(X_graphs_independent)):
    f = X_rho_formulas[j][0,1,2,3]
    v = vector(QQ, len(X_rho_basis), sparse=True)
    for coeff, monomial in zip(f.coefficients(), f.monomials()):
        monomial_index = X_rho_monomial_index[monomial]
        v[monomial_index] = coeff
    X_rho_evaluation_matrix.set_column(j, v)
print("Calculated X_rho_evaluation_matrix", flush=True)

print("Calculating rhodot_vector", flush=True)
rhodot_vector = vector(QQ, len(X_rho_basis), sparse=True)
for coeff, monomial in zip(rhodot.coefficients(), rhodot.monomials()):
    monomial_index = X_rho_monomial_index[monomial]
    rhodot_vector[monomial_index] = coeff
print("Calculated rhodot_vector", flush=True)

# Poisson-triviality $Q_{\text{tetra}}(P) = [\![P, X]\!]$ of the tetrahedral flow

# There is a linear combination of graphs $X$ such that $[\![P, X]\!]$ evaluates to $Q_{\text{tetra}}(P)$:

print("Calculating X_solution_vector", flush=True)
big_matrix = X_a_evaluation_matrix[0].stack(X_a_evaluation_matrix[1]).stack(X_rho_evaluation_matrix)
big_vector = vector(list(-adot_vector[0]) + list(-adot_vector[1]) + list(rhodot_vector))
X_solution_vector = big_matrix.solve_right(big_vector)
print("Calculated X_solution_vector", flush=True)
print("X_solution_vector =", X_solution_vector, flush=True)

X_solution = sum(c*f for c, f in zip(X_solution_vector, X_formulas_independent))

print("P.bracket(X_solution) == Q_tetra:", P.bracket(X_solution) == Q_tetra, flush=True)

# Parameters in the solution $X$ to $[\![P,X]\!] = Q_{\text{tetra}}(P)$

print("Number of parameters in the solution:", big_matrix.right_nullity(), flush=True)

print("Basis of kernel:", big_matrix.right_kernel().basis(), flush=True)
\end{verbatim}
\normalsize
The output:
\tiny
\begin{verbatim}
Number of graphs in X: 176
Calculating X_formulas
Calculated X_formulas
Number of monomials in components of X: [83946, 83946, 83946, 83946]
Calculating X_evaluation_matrix
Calculated X_evaluation_matrix
Number of linearly independent graphs in X: 86
Maximal subset of linearly independent graphs in X: [0, 1, 2, 3, 4, 5, 6, 7, 8, 9, 10, 11, 12, 13, 14, 15, 16, 17, 18, 19, 20, 21, 22, 23,
25, 26, 27, 28, 29, 30, 31, 32, 33, 34, 35, 36, 37, 40, 41, 42, 43, 44, 45, 47, 48, 49, 52, 53, 57, 78, 79, 80, 82, 83, 90, 91, 94, 96, 102,
103, 107, 129, 130, 131, 132, 133, 134, 135, 136, 138, 139, 140, 141, 142, 146, 147, 149, 150, 151, 158, 159, 161, 166, 168, 169, 170]
Calculating adot
Calculated adot
Calculating X_a_formulas
Calculated X_a_formulas
Number of monomials in X_a_basis: 103944 103944
Calculating X_a_evaluation_matrix
Calculated X_a_evaluation_matrix
Calculating adot_vector
Calculated adot_vector
Calculating Q_tetra
Calculated Q_tetra
Calculating rhodot
Calculated rhodot
Have nice expression for Q_tetra: True
Calculating X_rho_formulas
Calculated X_rho_formulas
Number of monomials in X_rho_basis: 334164
Calculating X_rho_evaluation_matrix
Calculated X_rho_evaluation_matrix
Calculating rhodot_vector
Calculated rhodot_vector
Calculating X_solution_vector
Traceback (most recent call last):
  File "/home2/s3058069/clean_4D_to_3D_vanishing_up.sage.py", line 314, in <module>
    X_solution_vector = big_matrix.solve_right(big_vector)
  File "sage/matrix/matrix2.pyx", line 939, in sage.matrix.matrix2.Matrix.solve_right (build/cythonized/sage/matrix/matrix2.c:16151)
  File "sage/matrix/matrix2.pyx", line 1062, in sage.matrix.matrix2.Matrix._solve_right_general (build/cythonized/sage/matrix/matrix2.c:17699)
ValueError: matrix equation has no solutions

###############################################################################
Hábrók Cluster

\end{verbatim}
\normalsize

\subsection{The 3D sunflower micro-graphs which yield linearly independent formulas do not extend to a 4D solution}\label{A5}

This result is in Idea \ref{ind}, Proposition \ref{fail4} and line 12 in Table~\ref{updowntable} on p.~\pageref{updowntable}

The script:
\tiny
\begin{verbatim}
# This is the "Plain 4D solution" code for 4D, meaning that we do NOT skew the encodings. Any solution given by this code is not a 'good' solution.
# Here, we take the encodings of the 20 linearly independent 3D sunflower graphs (linearly independent as formulas), and expand them to 4D
# We search over these for a 'not good' 4D solution

NPROCS=32

import warnings
warnings.filterwarnings("ignore", category=DeprecationWarning)

# Graphs for the vector field $X$

X_lin_ind_expanded_to_4D = [(0,1,4,7,1,3,5,8,1,2,6,9),(0,2,4,7,1,3,5,8,1,2,6,9)]
for i1 in [4,7]:
  for i2 in [4,7]:
    for j1 in [5,8]:
      for j2 in [5,8]:
        for k1 in [6,9]:
          for k2 in [6,9]:
            X_lin_ind_expanded_to_4D.append((0,1,4,7,1,k1,5,8,1,2,6,9))
            X_lin_ind_expanded_to_4D.append((0,1,4,7,1,3,5,8,i1,2,6,9))
            X_lin_ind_expanded_to_4D.append((0,1,4,7,1,k1,5,8,i1,2,6,9))
            X_lin_ind_expanded_to_4D.append((0,1,4,7,i1,k1,5,8,1,2,6,9))
            X_lin_ind_expanded_to_4D.append((0,1,4,7,i1,3,5,8,i2,2,6,9))
            X_lin_ind_expanded_to_4D.append((0,1,4,7,i1,k1,5,8,i2,2,6,9))
            X_lin_ind_expanded_to_4D.append((0,1,4,7,1,k1,5,8,1,j1,6,9))
            X_lin_ind_expanded_to_4D.append((0,1,4,7,i1,k1,5,8,i2,j1,6,9))
            X_lin_ind_expanded_to_4D.append((0,2,4,7,1,k1,5,8,1,2,6,9))
            X_lin_ind_expanded_to_4D.append((0,2,4,7,i1,3,5,8,i2,2,6,9))
            X_lin_ind_expanded_to_4D.append((0,2,4,7,i1,k1,5,8,i2,2,6,9))
            X_lin_ind_expanded_to_4D.append((0,2,4,7,1,3,5,8,1,j1,6,9))
            X_lin_ind_expanded_to_4D.append((0,2,4,7,1,k1,5,8,1,j1,6,9))
            X_lin_ind_expanded_to_4D.append((0,2,4,7,i1,3,5,8,1,j1,6,9))
            X_lin_ind_expanded_to_4D.append((0,j1,4,7,1,3,5,8,1,2,6,9))
            X_lin_ind_expanded_to_4D.append((0,j1,4,7,i1,3,5,8,i2,2,6,9))
            X_lin_ind_expanded_to_4D.append((0,j1,4,7,1,3,5,8,1,j2,6,9))
            X_lin_ind_expanded_to_4D.append((0,j1,4,7,1,k1,5,8,1,j2,6,9))

X_graph_encodings = list(set(X_lin_ind_expanded_to_4D))

# Convert the encodings to graphs:

from gcaops.graph.formality_graph import FormalityGraph
def encoding_to_graph(encoding):
    targets = [encoding[0:4], encoding[4:8], encoding[8:12]]
    edges = sum([[(k+1,v) for v in t] for (k,t) in enumerate(targets)], [])
    return FormalityGraph(1, 9, edges)

X_graphs = [encoding_to_graph(e) for e in X_graph_encodings]
print("Number of graphs in X:", len(X_graphs), flush=True)

# Formulas of graphs for $X$

# We want to evaluate each graph to a formula, inserting three copies of $\varrho\varepsilon^{ijk}$ 
# and three copies of the Casimir $a$ into the aerial vertices.

# First define the even coordinates $x,y,z$ (base), $\varrho, a$ (fibre) and odd coordinates $\xi_0, \xi_1, \xi_2$:

from gcaops.algebra.differential_polynomial_ring import DifferentialPolynomialRing
D4 = DifferentialPolynomialRing(QQ, ('rho','a1','a2'), ('x','y','z','w'), max_differential_orders=[3+1,1+3+1,1+3+1])
rho, a1, a2 = D4.fibre_variables()
x,y,z,w = D4.base_variables()
even_coords = [x,y,z,w]

from gcaops.algebra.superfunction_algebra import SuperfunctionAlgebra

S4.<xi0,xi1,xi2,xi3> = SuperfunctionAlgebra(D4, D4.base_variables())
xi = S4.gens()
odd_coords = xi
epsilon = xi[0]*xi[1]*xi[2]*xi[3] # Levi-Civita tensor

# Now evaluate each graph to a formula:

import itertools
from multiprocessing import Pool

def evaluate_graph(g):
    E = x*xi[0] + y*xi[1] + z*xi[2] + w*xi[3] # Euler vector field, to insert into ground vertex. Incoming derivative d/dx^i will result in xi[i].
    result = S4.zero()
    for index_choice in itertools.product(itertools.permutations(range(4)), repeat=3):
        sign = epsilon[index_choice[0]] * epsilon[index_choice[1]] * epsilon[index_choice[2]]
        # NOTE: This assumes the ground vertex is labeled 0, the vertices of out-degree 4 are labeled 1, 2, 3, 
        # and the Casimirs are labeled 4, 5, 6; 7, 8, 9.
        vertex_content = [E, S4(rho), S4(rho), S4(rho), S4(a1), S4(a1), S4(a1), S4(a2), S4(a2), S4(a2)]
        for ((source, target), index) in zip(g.edges(), sum(map(list, index_choice), [])):
            vertex_content[target] = vertex_content[target].derivative(even_coords[index])
        result += sign * prod(vertex_content)
    return result

print("Calculating X_formulas", flush=True)
X_formulas = []
with Pool(processes=NPROCS) as pool:
    X_formulas = list(pool.imap(evaluate_graph, X_graphs))
# NOTE: Necessary fixup of parents after multiprocessing.
X_formulas = [S4(X_formula) for X_formula in X_formulas]
print("Calculated X_formulas", flush=True)

# Each resulting formula is of the form $X_g = \sum_{i=1}^3 X_g^i \xi_i$ where $X_g^i$ is a differential polynomial depending on the graph $g$.

# Relations between (formulas of) graphs in $X$

# For each $i = 1,2,3,4$ we collect the distinct differential monomials in $X_g^i$ ranging over all $g$:

X_monomial_basis = [set([]) for i in range(4)]
for i in range(4):
    for X in X_formulas:
        X_monomial_basis[i] |= set(X[i].monomials())
X_monomial_basis = [list(b) for b in X_monomial_basis]
X_monomial_index = [{m : k for k, m in enumerate(b)} for b in X_monomial_basis]

print("Number of monomials in components of X:", [len(b) for b in X_monomial_basis], flush=True)

X_monomial_count = sum(len(b) for b in X_monomial_basis)

# Now the graph-to-formula evaluation can be expressed as a matrix:

X_evaluation_matrix = matrix(QQ, X_monomial_count, len(X_graphs), sparse=True)
for i in range(len(X_graphs)):
    v = vector(QQ, X_monomial_count, sparse=True)
    index_shift = 0
    for j in range(4):
        f = X_formulas[i][j]
        for coeff, monomial in zip(f.coefficients(), f.monomials()):
            monomial_index = X_monomial_index[j][monomial]
            v[index_shift + monomial_index] = coeff
        index_shift += len(X_monomial_basis[j])
    X_evaluation_matrix.set_column(i, v)

# Here is the number of linearly independent graphs:

print("Number of linearly independent graphs in X:", X_evaluation_matrix.rank(), flush=True)

# We collect the linearly independent graphs and their formulas in new shorter lists:

pivots = X_evaluation_matrix.pivots()
print("Maximal subset of linearly independent graphs:", list(pivots), flush=True)
X_graphs_independent = [X_graphs[k] for k in pivots]
X_formulas_independent = [X_formulas[k] for k in pivots]

# The tetrahedron operation $\operatorname{Op}(\gamma_3)$

from gcaops.graph.undirected_graph_complex import UndirectedGraphComplex

GC = UndirectedGraphComplex(QQ, implementation='vector', sparse=True)
tetrahedron = GC.cohomology_basis(4, 6)[0]

from gcaops.graph.directed_graph_complex import DirectedGraphComplex

dGC = DirectedGraphComplex(QQ, implementation='vector')
tetrahedron_oriented = dGC(tetrahedron)
tetrahedron_oriented_filtered = tetrahedron_oriented.filter(max_out_degree=2)
tetrahedron_operation = S4.graph_operation(tetrahedron_oriented_filtered)

# The flow $\dot{a}_i = 4\cdot\operatorname{Op}(\gamma_3)(P,P,P,a_i)$

# Let $P$ be the 4D Nambu–Poisson bracket with pre-factor $\varrho$ and Casimir functions $a_1, a_2$:

P = (rho*epsilon).bracket(a1).bracket(a2)

# We calculate $\dot{a}_i = 4\cdot\operatorname{Op}(\gamma_3)(P,P,P,a_i)$:

def casimir_flow(f):
    return 4*tetrahedron_operation(P,P,P,f)

a = [S4(a1), S4(a2)]

print("Calculating adot", flush=True)
adot = [casimir_flow(a_i) for a_i in a]
#adot = []
#with Pool(processes=2) as pool:
#    adot = list(pool.imap(casimir_flow, a))
# NOTE: Necessary fixup of parents after multiprocessing.
adot = [S4(D4(adot_i[()])) for adot_i in adot]
print("Calculated adot", flush=True)

# We calculate $X_g(a_i) = [\![X_g, a_i]\!]$ for $i=1,2$:

# TODO: Parallelize
print("Calculating X_a_formulas", flush=True)
X_a_formulas = [[X_formula.bracket(f) for X_formula in X_formulas_independent] for f in a]
print("Calculated X_a_formulas", flush=True)

# We express $g \mapsto X_g(a_i)$ as a sparse matrix for $i=1,2$:

X_a_basis = [set(f[()].monomials()) for f in adot]
for k in range(len(a)):
    for X_a_formula in X_a_formulas[k]:
        X_a_basis[k] |= set(X_a_formula[()].monomials())
X_a_basis = [list(B) for B in X_a_basis]
print("Number of monomials in X_a_basis:", len(X_a_basis[0]), len(X_a_basis[1]), flush=True)

print("Calculating X_a_evaluation_matrix", flush=True)
X_a_monomial_index = [{m : k for k, m in enumerate(B)} for B in X_a_basis]
X_a_evaluation_matrix = [matrix(QQ, len(B), len(X_graphs_independent), sparse=True) for B in X_a_basis]
for i in range(len(a)):
    for j in range(len(X_graphs_independent)):
        v = vector(QQ, len(X_a_basis[i]), sparse=True)
        f = X_a_formulas[i][j][()]
        for coeff, monomial in zip(f.coefficients(), f.monomials()):
            monomial_index = X_a_monomial_index[i][monomial]
            v[monomial_index] = coeff
        X_a_evaluation_matrix[i].set_column(j, v)
print("Calculated X_a_evaluation_matrix", flush=True)

# We express $\dot{a}_i$ as a vector for $i=1,2$:

print("Calculating adot_vector", flush=True)
adot_vector = [vector(QQ, len(B)) for B in X_a_basis]
for i in range(len(a)):
    f = adot[i][()]
    for coeff, monomial in zip(f.coefficients(), f.monomials()):
        monomial_index = X_a_monomial_index[i][monomial]
        adot_vector[i][monomial_index] = coeff
print("Calculated adot_vector", flush=True)

# The flow $\dot{\varrho}$

print("Calculating Q_tetra", flush=True)
Q_tetra = tetrahedron_operation(P,P,P,P)
print("Calculated Q_tetra", flush=True)

print("Calculating rhodot", flush=True)
P0 = (rho*epsilon).bracket(adot[0]).bracket(a2)
P1 = (rho*epsilon).bracket(a1).bracket(adot[1])
Q_remainder = Q_tetra - P0 - P1
P_withoutrho = epsilon.bracket(a1).bracket(a2)
rhodot = Q_remainder[0,1] // P_withoutrho[0,1]
print("Calculated rhodot", flush=True)

print("Have nice expression for Q_tetra:", Q_tetra == rhodot*P_withoutrho + P0 + P1, flush=True)

print("Calculating X_rho_formulas", flush=True)
# TODO: Parallelize
X_rho_formulas = [X_formula.bracket(rho*epsilon) for X_formula in X_formulas_independent]
print("Calculated X_rho_formulas", flush=True)

X_rho_basis = set(rhodot.monomials())
for X_rho_formula in X_rho_formulas:
    X_rho_basis |= set(X_rho_formula[0,1,2,3].monomials())
X_rho_basis = list(X_rho_basis)
print("Number of monomials in X_rho_basis:", len(X_rho_basis), flush=True)

print("Calculating X_rho_evaluation_matrix", flush=True)
X_rho_monomial_index = {m : k for k, m in enumerate(X_rho_basis)}
X_rho_evaluation_matrix = matrix(QQ, len(X_rho_basis), len(X_graphs_independent), sparse=True)
for j in range(len(X_graphs_independent)):
    f = X_rho_formulas[j][0,1,2,3]
    v = vector(QQ, len(X_rho_basis), sparse=True)
    for coeff, monomial in zip(f.coefficients(), f.monomials()):
        monomial_index = X_rho_monomial_index[monomial]
        v[monomial_index] = coeff
    X_rho_evaluation_matrix.set_column(j, v)
print("Calculated X_rho_evaluation_matrix", flush=True)

print("Calculating rhodot_vector", flush=True)
rhodot_vector = vector(QQ, len(X_rho_basis), sparse=True)
for coeff, monomial in zip(rhodot.coefficients(), rhodot.monomials()):
    monomial_index = X_rho_monomial_index[monomial]
    rhodot_vector[monomial_index] = coeff
print("Calculated rhodot_vector", flush=True)

# Poisson-triviality $Q_{\text{tetra}}(P) = [\![P, X]\!]$ of the tetrahedral flow

# There is a linear combination of graphs $X$ such that $[\![P, X]\!]$ evaluates to $Q_{\text{tetra}}(P)$:

print("Calculating X_solution_vector", flush=True)
big_matrix = X_a_evaluation_matrix[0].stack(X_a_evaluation_matrix[1]).stack(X_rho_evaluation_matrix)
big_vector = vector(list(-adot_vector[0]) + list(-adot_vector[1]) + list(-rhodot_vector))
X_solution_vector = big_matrix.solve_right(big_vector)
print("Calculated X_solution_vector", flush=True)
print("X_solution_vector =", X_solution_vector, flush=True)

X_solution = sum(c*f for c, f in zip(X_solution_vector, X_formulas_independent))

print("P.bracket(X_solution) == Q_tetra:", P.bracket(X_solution) == Q_tetra, flush=True)

# Parameters in the solution $X$ to $[\![P,X]\!] = Q_{\text{tetra}}(P)$

print("Number of parameters in the solution:", big_matrix.right_nullity(), flush=True)

print("Basis of kernel:", big_matrix.right_kernel().basis(), flush=True)    
\end{verbatim}
\normalsize
The output:
\tiny
\begin{verbatim}
Number of graphs in X: 92
Calculating X_formulas
Calculated X_formulas
Number of monomials in components of X: [82368, 82368, 82368, 82368]
Number of linearly independent graphs in X: 81
Maximal subset of linearly independent graphs: [0, 1, 2, 3, 4, 5, 6, 7, 8, 9, 10, 11, 13, 14, 15, 16, 17, 18, 19, 20,
21, 22, 23, 24, 25, 26, 27, 28, 30, 31, 32, 33, 34, 35, 36, 37, 38, 39, 41, 44, 45, 46, 47, 48, 49, 50, 51, 52, 53, 54,
55, 56, 57, 58, 59, 60, 61, 62, 63, 64, 65, 66, 67, 68, 70, 71, 72, 74, 75, 76, 78, 79, 80, 81, 84, 85, 86, 87, 89, 90, 91]
Calculating adot
Calculated adot
Calculating X_a_formulas
Calculated X_a_formulas
Number of monomials in X_a_basis: 101304 101304
Calculating X_a_evaluation_matrix
Calculated X_a_evaluation_matrix
Calculating adot_vector
Calculated adot_vector
Calculating Q_tetra
Calculated Q_tetra
Calculating rhodot
Calculated rhodot
Have nice expression for Q_tetra: True
Calculating X_rho_formulas
Calculated X_rho_formulas
Number of monomials in X_rho_basis: 344628
Calculating X_rho_evaluation_matrix
Calculated X_rho_evaluation_matrix
Calculating rhodot_vector
Calculated rhodot_vector
Calculating X_solution_vector
Traceback (most recent call last):
  File "/home2/s3058069/3D_lin_ind_expanded_to_4D_code.sage.py", line 271, in <module>
    X_solution_vector = big_matrix.solve_right(big_vector)
  File "sage/matrix/matrix2.pyx", line 939, in sage.matrix.matrix2.Matrix.solve_right (build/cythonized/sage/matrix/matrix2.c:16151)
  File "sage/matrix/matrix2.pyx", line 1062, in sage.matrix.matrix2.Matrix._solve_right_general (build/cythonized/sage/matrix/matrix2.c:17699)
ValueError: matrix equation has no solutions

###############################################################################
Hábrók Cluster
    
\end{verbatim}
\normalsize
\subsection{The 3D sunflower micro-graphs manually identified as extending to known 4D solution's components do extend to a 4D solution}\label{A6}

This result is in Ideas \ref{manual},\ref{success_idea}, Propositions \ref{17_3D},\ref{17_3D_give_4D_sol} and line 13 in Table~\ref{updowntable} on p.~\pageref{updowntable}.

The script:
\tiny
\begin{verbatim}
# This is the "Plain 4D solution" code for 4D, meaning that we do NOT skew the encodings. Any solution given by this code is not a 'good' solution,
# meaning it is not skew-symmetric with respect to swapping a1 and a2.
# Here, we take the encodings of the 17 3D sunflower graphs which give the 4D graphs, and expand them to 4D
# We search over these for a 'not good' 4D solution

NPROCS=32

import warnings
warnings.filterwarnings("ignore", category=DeprecationWarning)

# Graphs for the vector field $X$

X_lin_ind_expanded_to_4D = [(0,1,4,7,1,3,5,8,1,2,6,9),(0,2,4,7,1,3,5,8,1,2,6,9)]
for i1 in [4,7]:
  for i2 in [4,7]:
    for j1 in [5,8]:
      for j2 in [5,8]:
        for k1 in [6,9]:
          for k2 in [6,9]:
            X_lin_ind_expanded_to_4D.append((0,1,4,7,1,k1,5,8,i1,2,6,9))
            X_lin_ind_expanded_to_4D.append((0,1,4,7,i1,3,5,8,i2,2,6,9)) 
            X_lin_ind_expanded_to_4D.append((0,1,4,7,i1,k1,5,8,i2,2,6,9)) 
            X_lin_ind_expanded_to_4D.append((0,1,4,7,1,k1,5,8,1,j1,6,9)) 
            X_lin_ind_expanded_to_4D.append((0,1,4,7,i1,k1,5,8,i2,j1,6,9)) 
            X_lin_ind_expanded_to_4D.append((0,2,4,7,1,3,5,8,1,j1,6,9)) 
            X_lin_ind_expanded_to_4D.append((0,2,4,7,1,k1,5,8,1,j1,6,9)) 
            X_lin_ind_expanded_to_4D.append((0,2,4,7,i1,3,5,8,1,j1,6,9)) 
            X_lin_ind_expanded_to_4D.append((0,j1,4,7,1,3,5,8,1,2,6,9)) 
            X_lin_ind_expanded_to_4D.append((0,j1,4,7,1,k1,5,8,1,j2,6,9)) 
            X_lin_ind_expanded_to_4D.append((0,j1,4,7,1,k1,5,8,1,2,6,9)) # not lin ind and not vanishing
            X_lin_ind_expanded_to_4D.append((0,j1,4,7,i1,k1,5,8,i2,2,6,9)) # not lin ind and not vanishing
            X_lin_ind_expanded_to_4D.append((0,j1,4,7,1,k1,5,8,i2,2,6,9)) # vanishing	
            X_lin_ind_expanded_to_4D.append((0,2,4,7,1,k1,5,8,i2,j1,6,9)) # vanishing
            X_lin_ind_expanded_to_4D.append((0,j1,4,7,i1,k1,5,8,i2,j1,6,9)) # vanishing

		

X_graph_encodings = list(set(X_lin_ind_expanded_to_4D))


# Convert the encodings to graphs:

from gcaops.graph.formality_graph import FormalityGraph
def encoding_to_graph(encoding):
    targets = [encoding[0:4], encoding[4:8], encoding[8:12]]
    edges = sum([[(k+1,v) for v in t] for (k,t) in enumerate(targets)], [])
    return FormalityGraph(1, 9, edges)

X_graphs = [encoding_to_graph(e) for e in X_graph_encodings]
print("Number of graphs in X:", len(X_graphs), flush=True)

# Formulas of graphs for $X$

# We want to evaluate each graph to a formula, inserting three copies of $\varrho\varepsilon^{ijk}$ 
# and three copies of the Casimir $a$ into the aerial vertices.

# First define the even coordinates $x,y,z$ (base), $\varrho, a$ (fibre) and odd coordinates $\xi_0, \xi_1, \xi_2$:

from gcaops.algebra.differential_polynomial_ring import DifferentialPolynomialRing
D4 = DifferentialPolynomialRing(QQ, ('rho','a1','a2'), ('x','y','z','w'), max_differential_orders=[3+1,1+3+1,1+3+1])
rho, a1, a2 = D4.fibre_variables()
x,y,z,w = D4.base_variables()
even_coords = [x,y,z,w]

from gcaops.algebra.superfunction_algebra import SuperfunctionAlgebra

S4.<xi0,xi1,xi2,xi3> = SuperfunctionAlgebra(D4, D4.base_variables())
xi = S4.gens()
odd_coords = xi
epsilon = xi[0]*xi[1]*xi[2]*xi[3] # Levi-Civita tensor

# Now evaluate each graph to a formula:

import itertools
from multiprocessing import Pool

def evaluate_graph(g):
    E = x*xi[0] + y*xi[1] + z*xi[2] + w*xi[3] # Euler vector field, to insert into ground vertex. Incoming derivative d/dx^i will result in xi[i].
    result = S4.zero()
    for index_choice in itertools.product(itertools.permutations(range(4)), repeat=3):
        sign = epsilon[index_choice[0]] * epsilon[index_choice[1]] * epsilon[index_choice[2]]
        # NOTE: This assumes the ground vertex is labeled 0, the vertices of out-degree 4 are labeled 1, 2, 3, 
        # and the Casimirs are labeled 4, 5, 6; 7, 8, 9.
        vertex_content = [E, S4(rho), S4(rho), S4(rho), S4(a1), S4(a1), S4(a1), S4(a2), S4(a2), S4(a2)]
        for ((source, target), index) in zip(g.edges(), sum(map(list, index_choice), [])):
            vertex_content[target] = vertex_content[target].derivative(even_coords[index])
        result += sign * prod(vertex_content)
    return result

print("Calculating X_formulas", flush=True)
X_formulas = []
with Pool(processes=NPROCS) as pool:
    X_formulas = list(pool.imap(evaluate_graph, X_graphs))
# NOTE: Necessary fixup of parents after multiprocessing.
X_formulas = [S4(X_formula) for X_formula in X_formulas]
print("Calculated X_formulas", flush=True)

# Each resulting formula is of the form $X_g = \sum_{i=1}^3 X_g^i \xi_i$ where $X_g^i$ is a differential polynomial depending on the graph $g$.

# Relations between (formulas of) graphs in $X$

# For each $i = 1,2,3,4$ we collect the distinct differential monomials in $X_g^i$ ranging over all $g$:

X_monomial_basis = [set([]) for i in range(4)]
for i in range(4):
    for X in X_formulas:
        X_monomial_basis[i] |= set(X[i].monomials())
X_monomial_basis = [list(b) for b in X_monomial_basis]
X_monomial_index = [{m : k for k, m in enumerate(b)} for b in X_monomial_basis]

print("Number of monomials in components of X:", [len(b) for b in X_monomial_basis], flush=True)

X_monomial_count = sum(len(b) for b in X_monomial_basis)

# Now the graph-to-formula evaluation can be expressed as a matrix:

X_evaluation_matrix = matrix(QQ, X_monomial_count, len(X_graphs), sparse=True)
for i in range(len(X_graphs)):
    v = vector(QQ, X_monomial_count, sparse=True)
    index_shift = 0
    for j in range(4):
        f = X_formulas[i][j]
        for coeff, monomial in zip(f.coefficients(), f.monomials()):
            monomial_index = X_monomial_index[j][monomial]
            v[index_shift + monomial_index] = coeff
        index_shift += len(X_monomial_basis[j])
    X_evaluation_matrix.set_column(i, v)

# Here is the number of linearly independent graphs:

print("Number of linearly independent graphs in X:", X_evaluation_matrix.rank(), flush=True)

# We collect the linearly independent graphs and their formulas in new shorter lists:

pivots = X_evaluation_matrix.pivots()
print("Maximal subset of linearly independent graphs:", list(pivots), flush=True)
X_graphs_independent = [X_graphs[k] for k in pivots]
X_formulas_independent = [X_formulas[k] for k in pivots]

# The tetrahedron operation $\operatorname{Op}(\gamma_3)$

from gcaops.graph.undirected_graph_complex import UndirectedGraphComplex

GC = UndirectedGraphComplex(QQ, implementation='vector', sparse=True)
tetrahedron = GC.cohomology_basis(4, 6)[0]

from gcaops.graph.directed_graph_complex import DirectedGraphComplex

dGC = DirectedGraphComplex(QQ, implementation='vector')
tetrahedron_oriented = dGC(tetrahedron)
tetrahedron_oriented_filtered = tetrahedron_oriented.filter(max_out_degree=2)
tetrahedron_operation = S4.graph_operation(tetrahedron_oriented_filtered)

# The flow $\dot{a}_i = 4\cdot\operatorname{Op}(\gamma_3)(P,P,P,a_i)$

# Let $P$ be the 4D Nambu–Poisson bracket with pre-factor $\varrho$ and Casimir functions $a_1, a_2$:

P = (rho*epsilon).bracket(a1).bracket(a2)

# We calculate $\dot{a}_i = 4\cdot\operatorname{Op}(\gamma_3)(P,P,P,a_i)$:

def casimir_flow(f):
    return 4*tetrahedron_operation(P,P,P,f)

a = [S4(a1), S4(a2)]

print("Calculating adot", flush=True)
adot = [casimir_flow(a_i) for a_i in a]
#adot = []
#with Pool(processes=2) as pool:
#    adot = list(pool.imap(casimir_flow, a))
# NOTE: Necessary fixup of parents after multiprocessing.
adot = [S4(D4(adot_i[()])) for adot_i in adot]
print("Calculated adot", flush=True)

# We calculate $X_g(a_i) = [\![X_g, a_i]\!]$ for $i=1,2$:

# TODO: Parallelize
print("Calculating X_a_formulas", flush=True)
X_a_formulas = [[X_formula.bracket(f) for X_formula in X_formulas_independent] for f in a]
print("Calculated X_a_formulas", flush=True)

# We express $g \mapsto X_g(a_i)$ as a sparse matrix for $i=1,2$:

X_a_basis = [set(f[()].monomials()) for f in adot]
for k in range(len(a)):
    for X_a_formula in X_a_formulas[k]:
        X_a_basis[k] |= set(X_a_formula[()].monomials())
X_a_basis = [list(B) for B in X_a_basis]
print("Number of monomials in X_a_basis:", len(X_a_basis[0]), len(X_a_basis[1]), flush=True)

print("Calculating X_a_evaluation_matrix", flush=True)
X_a_monomial_index = [{m : k for k, m in enumerate(B)} for B in X_a_basis]
X_a_evaluation_matrix = [matrix(QQ, len(B), len(X_graphs_independent), sparse=True) for B in X_a_basis]
for i in range(len(a)):
    for j in range(len(X_graphs_independent)):
        v = vector(QQ, len(X_a_basis[i]), sparse=True)
        f = X_a_formulas[i][j][()]
        for coeff, monomial in zip(f.coefficients(), f.monomials()):
            monomial_index = X_a_monomial_index[i][monomial]
            v[monomial_index] = coeff
        X_a_evaluation_matrix[i].set_column(j, v)
print("Calculated X_a_evaluation_matrix", flush=True)

# We express $\dot{a}_i$ as a vector for $i=1,2$:

print("Calculating adot_vector", flush=True)
adot_vector = [vector(QQ, len(B)) for B in X_a_basis]
for i in range(len(a)):
    f = adot[i][()]
    for coeff, monomial in zip(f.coefficients(), f.monomials()):
        monomial_index = X_a_monomial_index[i][monomial]
        adot_vector[i][monomial_index] = coeff
print("Calculated adot_vector", flush=True)

# The flow $\dot{\varrho}$

print("Calculating Q_tetra", flush=True)
Q_tetra = tetrahedron_operation(P,P,P,P)
print("Calculated Q_tetra", flush=True)

print("Calculating rhodot", flush=True)
P0 = (rho*epsilon).bracket(adot[0]).bracket(a2)
P1 = (rho*epsilon).bracket(a1).bracket(adot[1])
Q_remainder = Q_tetra - P0 - P1
P_withoutrho = epsilon.bracket(a1).bracket(a2)
rhodot = Q_remainder[0,1] // P_withoutrho[0,1]
print("Calculated rhodot", flush=True)

print("Have nice expression for Q_tetra:", Q_tetra == rhodot*P_withoutrho + P0 + P1, flush=True)

print("Calculating X_rho_formulas", flush=True)
# TODO: Parallelize
X_rho_formulas = [X_formula.bracket(rho*epsilon) for X_formula in X_formulas_independent]
print("Calculated X_rho_formulas", flush=True)

X_rho_basis = set(rhodot.monomials())
for X_rho_formula in X_rho_formulas:
    X_rho_basis |= set(X_rho_formula[0,1,2,3].monomials())
X_rho_basis = list(X_rho_basis)
print("Number of monomials in X_rho_basis:", len(X_rho_basis), flush=True)

print("Calculating X_rho_evaluation_matrix", flush=True)
X_rho_monomial_index = {m : k for k, m in enumerate(X_rho_basis)}
X_rho_evaluation_matrix = matrix(QQ, len(X_rho_basis), len(X_graphs_independent), sparse=True)
for j in range(len(X_graphs_independent)):
    f = X_rho_formulas[j][0,1,2,3]
    v = vector(QQ, len(X_rho_basis), sparse=True)
    for coeff, monomial in zip(f.coefficients(), f.monomials()):
        monomial_index = X_rho_monomial_index[monomial]
        v[monomial_index] = coeff
    X_rho_evaluation_matrix.set_column(j, v)
print("Calculated X_rho_evaluation_matrix", flush=True)

print("Calculating rhodot_vector", flush=True)
rhodot_vector = vector(QQ, len(X_rho_basis), sparse=True)
for coeff, monomial in zip(rhodot.coefficients(), rhodot.monomials()):
    monomial_index = X_rho_monomial_index[monomial]
    rhodot_vector[monomial_index] = coeff
print("Calculated rhodot_vector", flush=True)

# Poisson-triviality $Q_{\text{tetra}}(P) = [\![P, X]\!]$ of the tetrahedral flow

# There is a linear combination of graphs $X$ such that $[\![P, X]\!]$ evaluates to $Q_{\text{tetra}}(P)$:

print("Calculating X_solution_vector", flush=True)
big_matrix = X_a_evaluation_matrix[0].stack(X_a_evaluation_matrix[1]).stack(X_rho_evaluation_matrix)
big_vector = vector(list(-adot_vector[0]) + list(-adot_vector[1]) + list(-rhodot_vector))
X_solution_vector = big_matrix.solve_right(big_vector)
print("Calculated X_solution_vector", flush=True)
print("X_solution_vector =", X_solution_vector, flush=True)

X_solution = sum(c*f for c, f in zip(X_solution_vector, X_formulas_independent))

print("P.bracket(X_solution) == Q_tetra:", P.bracket(X_solution) == Q_tetra, flush=True)

# Parameters in the solution $X$ to $[\![P,X]\!] = Q_{\text{tetra}}(P)$

print("Number of parameters in the solution:", big_matrix.right_nullity(), flush=True)

print("Basis of kernel:", big_matrix.right_kernel().basis(), flush=True)
\end{verbatim}
\normalsize
The output:
\tiny
\begin{verbatim}
Number of graphs in X: 110
Calculating X_formulas
Calculated X_formulas
Number of monomials in components of X: [81402, 81402, 81402, 81402]
Number of linearly independent graphs in X: 76
Maximal subset of linearly independent graphs: [0, 1, 2, 3, 4, 5, 6, 9, 10, 11, 12, 13, 14, 15, 16, 17, 20, 21, 22, 23, 24, 25,
27, 28, 29, 30, 31, 32, 33, 34, 35, 37, 39, 40, 42, 43, 45, 46, 48, 53, 55, 57, 58, 59, 60, 61, 62, 63, 64, 65, 66, 67, 69, 70,
71, 73, 74, 75, 76, 77, 78, 79, 82, 84, 85, 88, 90, 91, 92, 93, 96, 98, 101, 102, 105, 109]
Calculating adot
Calculated adot
Calculating X_a_formulas
Calculated X_a_formulas
Number of monomials in X_a_basis: 79908 79908
Calculating X_a_evaluation_matrix
Calculated X_a_evaluation_matrix
Calculating adot_vector
Calculated adot_vector
Calculating Q_tetra
Calculated Q_tetra
Calculating rhodot
Calculated rhodot
Have nice expression for Q_tetra: True
Calculating X_rho_formulas
Calculated X_rho_formulas
Number of monomials in X_rho_basis: 335028
Calculating X_rho_evaluation_matrix
Calculated X_rho_evaluation_matrix
Calculating rhodot_vector
Calculated rhodot_vector
Calculating X_solution_vector
Calculated X_solution_vector
X_solution_vector = (-24, 24, 24, -12, 24, 0, 0, -24, -24, 0, 0, -24, 0, 0, 24, -24, -72, -24, -24, -24, 0, 12, -24, -48, -24, 0,
6, -24, 0, 0, -16, -24, 24, 0, -16, 24, -16, 12, -24, 0, 12, 0, -24, -12, 6, 0, -16, 0, -8, 0, -48, -16, 0, -24, -8, -12, 0, 0, 24, 0,
24, -24, 0, -24, -24, -24, 0, 0, -8, 0, -16, -16, -16, -24, 0, -24)
P.bracket(X_solution) == Q_tetra: True
Number of parameters in the solution: 3
Basis of kernel: [
(0, 0, 0, 0, 0, 0, 0, 0, 0, 0, 0, 0, 0, 0, 0, 0, 0, 0, 0, 0, 0, 1, 0, 0, 0, 0, 0, 0, 0, 0, 0, 0, 1, 0, 0, 0, 0, 0, 0, 0, 1/2, 1/2,
0, 0, 0, 0, 0, 0, 0, 0, 0, 0, 0, 0, 0, 0, 0, 0, 0, 0, 0, 0, 0, 0, 0, 0, 0, 0, 0, 0, 0, 0, 0, 0, 0, 0),
(0, 0, 0, 0, 0, 0, 0, 0, 0, 0, 0, 0, 0, 0, 0, 0, 0, 0, 0, 0, 0, 0, 0, 0, 0, 0, 1, 0, 0, 0, 0, 0, 0, 0, 0, 0, 0, 2, 0, 0, 0, 0, 0, 0, 1,
0, 0, 0, 0, 2, 0, 0, 0, 0, 0, 0, 0, 0, 0, 0, 0, 0, 0, 0, 0, 0, 0, 0, 0, 0, 0, 0, 0, 0, 0, 0),
(0, 0, 0, 0, 0, 0, 0, 0, 0, 0, 0, 0, 0, 0, 0, 0, 0, 0, 0, 0, 0, 0, 0, 0, 0, 0, 0, 0, 0, 1, 0, 0, 0, 0, 0, 0, 0, 0, 0, 0, 0, 0,
0, 0, 0, 0, 0, 0, 0, 0, 0, 0, 0, 0, 0, 0, 1, 0, 0, 1, 0, 0, 1, 0, 0, 0, 0, 1, 0, 0, 0, 0, 0, 0, 1, 0)
]

###############################################################################
Hábrók Cluster

\end{verbatim}
\normalsize
\subsection{The 3D basic sunflower micro-graphs and the 3D vanishing sunflower micro-graphs extend to a 4D solution}\label{A7}

This result is in Idea \ref{general}, Proposition \ref{lin_ind_and_vanishing_to_4d} and line 14 in Table~\ref{updowntable} on p.~\pageref{updowntable}.

The script:
\tiny
\begin{verbatim}
# This is the "Plain 4D solution" code for 4D, meaning that we do NOT skew the encodings. Any solution given by this code is not a 'good' solution.
# Here, we take the encodings of the 20 linearly independent 3D sunflower graphs (linearly independent as formulas), and expand them to 4D
# To these expanded encodings we add the 3D vanishing encodings expanded to 4D
# We search over these for a 'not good' 4D solution

NPROCS=32

import warnings
warnings.filterwarnings("ignore", category=DeprecationWarning)

# Graphs for the vector field $X$

X_lin_ind_expanded_to_4D = [(0,1,4,7,1,3,5,8,1,2,6,9),(0,2,4,7,1,3,5,8,1,2,6,9)]
for i1 in [4,7]:
  for i2 in [4,7]:
    for j1 in [5,8]:
      for j2 in [5,8]:
        for k1 in [6,9]:
          for k2 in [6,9]:
            X_lin_ind_expanded_to_4D.append((0,1,4,7,1,k1,5,8,1,2,6,9))
            X_lin_ind_expanded_to_4D.append((0,1,4,7,1,3,5,8,i1,2,6,9))
            X_lin_ind_expanded_to_4D.append((0,1,4,7,1,k1,5,8,i1,2,6,9))
            X_lin_ind_expanded_to_4D.append((0,1,4,7,i1,k1,5,8,1,2,6,9))
            X_lin_ind_expanded_to_4D.append((0,1,4,7,i1,3,5,8,i2,2,6,9))
            X_lin_ind_expanded_to_4D.append((0,1,4,7,i1,k1,5,8,i2,2,6,9))
            X_lin_ind_expanded_to_4D.append((0,1,4,7,1,k1,5,8,1,j1,6,9))
            X_lin_ind_expanded_to_4D.append((0,1,4,7,i1,k1,5,8,i2,j1,6,9))
            X_lin_ind_expanded_to_4D.append((0,2,4,7,1,k1,5,8,1,2,6,9))
            X_lin_ind_expanded_to_4D.append((0,2,4,7,i1,3,5,8,i2,2,6,9))
            X_lin_ind_expanded_to_4D.append((0,2,4,7,i1,k1,5,8,i2,2,6,9))
            X_lin_ind_expanded_to_4D.append((0,2,4,7,1,3,5,8,1,j1,6,9))
            X_lin_ind_expanded_to_4D.append((0,2,4,7,1,k1,5,8,1,j1,6,9))
            X_lin_ind_expanded_to_4D.append((0,2,4,7,i1,3,5,8,1,j1,6,9))
            X_lin_ind_expanded_to_4D.append((0,j1,4,7,1,3,5,8,1,2,6,9))
            X_lin_ind_expanded_to_4D.append((0,j1,4,7,i1,3,5,8,i2,2,6,9))
            X_lin_ind_expanded_to_4D.append((0,j1,4,7,1,3,5,8,1,j2,6,9))
            X_lin_ind_expanded_to_4D.append((0,j1,4,7,1,k1,5,8,1,j2,6,9))

X_graph_encodings = list(set(X_lin_ind_expanded_to_4D))

vanishing_expanded=[(0,2,4,7,1,3,5,8,4,2,6,9),(0,2,4,7,1,3,5,8,7,2,6,9),(0,5,4,7,4,3,5,8,1,2,6,9),(0,8,4,7,4,3,5,8,1,2,6,9),
(0,8,4,7,7,3,5,8,1,2,6,9),(0,5,4,7,7,3,5,8,1,2,6,9),(0,5,4,7,1,3,5,8,4,2,6,9),(0,8,4,7,1,3,5,8,4,2,6,9),(0,8,4,7,1,3,5,8,7,2,6,9),
(0,5,4,7,1,3,5,8,7,2,6,9),(0,2,4,7,1,3,5,8,4,5,6,9),(0,2,4,7,1,3,5,8,7,5,6,9),(0,2,4,7,1,3,5,8,7,8,6,9),(0,2,4,7,1,3,5,8,4,8,6,9)]

for i in [4,7]:
    for j in [5,8]:
        for k in [6,9]:
            vanishing_expanded.append((0,1,4,7,1,k,5,8,i,j,6,9))
            vanishing_expanded.append((0,2,4,7,1,k,5,8,i,j,6,9))
            vanishing_expanded.append((0,2,4,7,i,k,5,8,1,j,6,9))
            vanishing_expanded.append((0,j,4,7,i,k,5,8,1,2,6,9))
            vanishing_expanded.append((0,j,4,7,1,k,5,8,i,2,6,9))
            vanishing_expanded.append((0,1,4,7,i,k,5,8,1,j,6,9))
            
for i in [4,7]:
    for j1 in [5,8]:
        for j2 in [5,8]:
            vanishing_expanded.append((0,j1,4,7,1,3,5,8,i,j2,6,9))
            
for i1 in [4,7]:
    for i2 in [4,7]:
        for j1 in [5,8]:
            for j2 in [5,8]:
                vanishing_expanded.append((0,j1,4,7,i1,3,5,8,i2,j2,6,9))
                
for i1 in [4,7]:
    for i2 in [4,7]:
        for j1 in [5,8]:
            for j2 in [5,8]:
                for k in [6,9]:
                    vanishing_expanded.append((0,j1,4,7,i1,k,5,8,i2,j2,6,9))

X_graph_encodings += vanishing_expanded


# Convert the encodings to graphs:

from gcaops.graph.formality_graph import FormalityGraph
def encoding_to_graph(encoding):
    targets = [encoding[0:4], encoding[4:8], encoding[8:12]]
    edges = sum([[(k+1,v) for v in t] for (k,t) in enumerate(targets)], [])
    return FormalityGraph(1, 9, edges)

X_graphs = [encoding_to_graph(e) for e in X_graph_encodings]
print("Number of graphs in X:", len(X_graphs), flush=True)

# Formulas of graphs for $X$

# We want to evaluate each graph to a formula, inserting three copies of $\varrho\varepsilon^{ijk}$ 
# and three copies of the Casimir $a$ into the aerial vertices.

# First define the even coordinates $x,y,z$ (base), $\varrho, a$ (fibre) and odd coordinates $\xi_0, \xi_1, \xi_2$:

from gcaops.algebra.differential_polynomial_ring import DifferentialPolynomialRing
D4 = DifferentialPolynomialRing(QQ, ('rho','a1','a2'), ('x','y','z','w'), max_differential_orders=[3+1,1+3+1,1+3+1])
rho, a1, a2 = D4.fibre_variables()
x,y,z,w = D4.base_variables()
even_coords = [x,y,z,w]

from gcaops.algebra.superfunction_algebra import SuperfunctionAlgebra

S4.<xi0,xi1,xi2,xi3> = SuperfunctionAlgebra(D4, D4.base_variables())
xi = S4.gens()
odd_coords = xi
epsilon = xi[0]*xi[1]*xi[2]*xi[3] # Levi-Civita tensor

# Now evaluate each graph to a formula:

import itertools
from multiprocessing import Pool

def evaluate_graph(g):
    E = x*xi[0] + y*xi[1] + z*xi[2] + w*xi[3] # Euler vector field, to insert into ground vertex. Incoming derivative d/dx^i will result in xi[i].
    result = S4.zero()
    for index_choice in itertools.product(itertools.permutations(range(4)), repeat=3):
        sign = epsilon[index_choice[0]] * epsilon[index_choice[1]] * epsilon[index_choice[2]]
        # NOTE: This assumes the ground vertex is labeled 0, the vertices of out-degree 4 are labeled 1, 2, 3, 
        # and the Casimirs are labeled 4, 5, 6; 7, 8, 9.
        vertex_content = [E, S4(rho), S4(rho), S4(rho), S4(a1), S4(a1), S4(a1), S4(a2), S4(a2), S4(a2)]
        for ((source, target), index) in zip(g.edges(), sum(map(list, index_choice), [])):
            vertex_content[target] = vertex_content[target].derivative(even_coords[index])
        result += sign * prod(vertex_content)
    return result

print("Calculating X_formulas", flush=True)
X_formulas = []
with Pool(processes=NPROCS) as pool:
    X_formulas = list(pool.imap(evaluate_graph, X_graphs))
# NOTE: Necessary fixup of parents after multiprocessing.
X_formulas = [S4(X_formula) for X_formula in X_formulas]
print("Calculated X_formulas", flush=True)

# Each resulting formula is of the form $X_g = \sum_{i=1}^3 X_g^i \xi_i$ where $X_g^i$ is a differential polynomial depending on the graph $g$.

# Relations between (formulas of) graphs in $X$

# For each $i = 1,2,3,4$ we collect the distinct differential monomials in $X_g^i$ ranging over all $g$:

X_monomial_basis = [set([]) for i in range(4)]
for i in range(4):
    for X in X_formulas:
        X_monomial_basis[i] |= set(X[i].monomials())
X_monomial_basis = [list(b) for b in X_monomial_basis]
X_monomial_index = [{m : k for k, m in enumerate(b)} for b in X_monomial_basis]

print("Number of monomials in components of X:", [len(b) for b in X_monomial_basis], flush=True)

X_monomial_count = sum(len(b) for b in X_monomial_basis)

# Now the graph-to-formula evaluation can be expressed as a matrix:

X_evaluation_matrix = matrix(QQ, X_monomial_count, len(X_graphs), sparse=True)
for i in range(len(X_graphs)):
    v = vector(QQ, X_monomial_count, sparse=True)
    index_shift = 0
    for j in range(4):
        f = X_formulas[i][j]
        for coeff, monomial in zip(f.coefficients(), f.monomials()):
            monomial_index = X_monomial_index[j][monomial]
            v[index_shift + monomial_index] = coeff
        index_shift += len(X_monomial_basis[j])
    X_evaluation_matrix.set_column(i, v)

# Here is the number of linearly independent graphs:

print("Number of linearly independent graphs in X:", X_evaluation_matrix.rank(), flush=True)

# We collect the linearly independent graphs and their formulas in new shorter lists:

pivots = X_evaluation_matrix.pivots()
print("Maximal subset of linearly independent graphs:", list(pivots), flush=True)
X_graphs_independent = [X_graphs[k] for k in pivots]
X_formulas_independent = [X_formulas[k] for k in pivots]

# The tetrahedron operation $\operatorname{Op}(\gamma_3)$

from gcaops.graph.undirected_graph_complex import UndirectedGraphComplex

GC = UndirectedGraphComplex(QQ, implementation='vector', sparse=True)
tetrahedron = GC.cohomology_basis(4, 6)[0]

from gcaops.graph.directed_graph_complex import DirectedGraphComplex

dGC = DirectedGraphComplex(QQ, implementation='vector')
tetrahedron_oriented = dGC(tetrahedron)
tetrahedron_oriented_filtered = tetrahedron_oriented.filter(max_out_degree=2)
tetrahedron_operation = S4.graph_operation(tetrahedron_oriented_filtered)

# The flow $\dot{a}_i = 4\cdot\operatorname{Op}(\gamma_3)(P,P,P,a_i)$

# Let $P$ be the 4D Nambu–Poisson bracket with pre-factor $\varrho$ and Casimir functions $a_1, a_2$:

P = (rho*epsilon).bracket(a1).bracket(a2)

# We calculate $\dot{a}_i = 4\cdot\operatorname{Op}(\gamma_3)(P,P,P,a_i)$:

def casimir_flow(f):
    return 4*tetrahedron_operation(P,P,P,f)

a = [S4(a1), S4(a2)]

print("Calculating adot", flush=True)
adot = [casimir_flow(a_i) for a_i in a]
#adot = []
#with Pool(processes=2) as pool:
#    adot = list(pool.imap(casimir_flow, a))
# NOTE: Necessary fixup of parents after multiprocessing.
adot = [S4(D4(adot_i[()])) for adot_i in adot]
print("Calculated adot", flush=True)

# We calculate $X_g(a_i) = [\![X_g, a_i]\!]$ for $i=1,2$:

# TODO: Parallelize
print("Calculating X_a_formulas", flush=True)
X_a_formulas = [[X_formula.bracket(f) for X_formula in X_formulas_independent] for f in a]
print("Calculated X_a_formulas", flush=True)

# We express $g \mapsto X_g(a_i)$ as a sparse matrix for $i=1,2$:

X_a_basis = [set(f[()].monomials()) for f in adot]
for k in range(len(a)):
    for X_a_formula in X_a_formulas[k]:
        X_a_basis[k] |= set(X_a_formula[()].monomials())
X_a_basis = [list(B) for B in X_a_basis]
print("Number of monomials in X_a_basis:", len(X_a_basis[0]), len(X_a_basis[1]), flush=True)

print("Calculating X_a_evaluation_matrix", flush=True)
X_a_monomial_index = [{m : k for k, m in enumerate(B)} for B in X_a_basis]
X_a_evaluation_matrix = [matrix(QQ, len(B), len(X_graphs_independent), sparse=True) for B in X_a_basis]
for i in range(len(a)):
    for j in range(len(X_graphs_independent)):
        v = vector(QQ, len(X_a_basis[i]), sparse=True)
        f = X_a_formulas[i][j][()]
        for coeff, monomial in zip(f.coefficients(), f.monomials()):
            monomial_index = X_a_monomial_index[i][monomial]
            v[monomial_index] = coeff
        X_a_evaluation_matrix[i].set_column(j, v)
print("Calculated X_a_evaluation_matrix", flush=True)

# We express $\dot{a}_i$ as a vector for $i=1,2$:

print("Calculating adot_vector", flush=True)
adot_vector = [vector(QQ, len(B)) for B in X_a_basis]
for i in range(len(a)):
    f = adot[i][()]
    for coeff, monomial in zip(f.coefficients(), f.monomials()):
        monomial_index = X_a_monomial_index[i][monomial]
        adot_vector[i][monomial_index] = coeff
print("Calculated adot_vector", flush=True)

# The flow $\dot{\varrho}$

print("Calculating Q_tetra", flush=True)
Q_tetra = tetrahedron_operation(P,P,P,P)
print("Calculated Q_tetra", flush=True)

print("Calculating rhodot", flush=True)
P0 = (rho*epsilon).bracket(adot[0]).bracket(a2)
P1 = (rho*epsilon).bracket(a1).bracket(adot[1])
Q_remainder = Q_tetra - P0 - P1
P_withoutrho = epsilon.bracket(a1).bracket(a2)
rhodot = Q_remainder[0,1] // P_withoutrho[0,1]
print("Calculated rhodot", flush=True)

print("Have nice expression for Q_tetra:", Q_tetra == rhodot*P_withoutrho + P0 + P1, flush=True)

print("Calculating X_rho_formulas", flush=True)
# TODO: Parallelize
X_rho_formulas = [X_formula.bracket(rho*epsilon) for X_formula in X_formulas_independent]
print("Calculated X_rho_formulas", flush=True)

X_rho_basis = set(rhodot.monomials())
for X_rho_formula in X_rho_formulas:
    X_rho_basis |= set(X_rho_formula[0,1,2,3].monomials())
X_rho_basis = list(X_rho_basis)
print("Number of monomials in X_rho_basis:", len(X_rho_basis), flush=True)

print("Calculating X_rho_evaluation_matrix", flush=True)
X_rho_monomial_index = {m : k for k, m in enumerate(X_rho_basis)}
X_rho_evaluation_matrix = matrix(QQ, len(X_rho_basis), len(X_graphs_independent), sparse=True)
for j in range(len(X_graphs_independent)):
    f = X_rho_formulas[j][0,1,2,3]
    v = vector(QQ, len(X_rho_basis), sparse=True)
    for coeff, monomial in zip(f.coefficients(), f.monomials()):
        monomial_index = X_rho_monomial_index[monomial]
        v[monomial_index] = coeff
    X_rho_evaluation_matrix.set_column(j, v)
print("Calculated X_rho_evaluation_matrix", flush=True)

print("Calculating rhodot_vector", flush=True)
rhodot_vector = vector(QQ, len(X_rho_basis), sparse=True)
for coeff, monomial in zip(rhodot.coefficients(), rhodot.monomials()):
    monomial_index = X_rho_monomial_index[monomial]
    rhodot_vector[monomial_index] = coeff
print("Calculated rhodot_vector", flush=True)

# Poisson-triviality $Q_{\text{tetra}}(P) = [\![P, X]\!]$ of the tetrahedral flow

# There is a linear combination of graphs $X$ such that $[\![P, X]\!]$ evaluates to $Q_{\text{tetra}}(P)$:

print("Calculating X_solution_vector", flush=True)
big_matrix = X_a_evaluation_matrix[0].stack(X_a_evaluation_matrix[1]).stack(X_rho_evaluation_matrix)
big_vector = vector(list(-adot_vector[0]) + list(-adot_vector[1]) + list(-rhodot_vector))
X_solution_vector = big_matrix.solve_right(big_vector)
print("Calculated X_solution_vector", flush=True)
print("X_solution_vector =", X_solution_vector, flush=True)

X_solution = sum(c*f for c, f in zip(X_solution_vector, X_formulas_independent))

print("P.bracket(X_solution) == Q_tetra:", P.bracket(X_solution) == Q_tetra, flush=True)

# Parameters in the solution $X$ to $[\![P,X]\!] = Q_{\text{tetra}}(P)$

print("Number of parameters in the solution:", big_matrix.right_nullity(), flush=True)

print("Basis of kernel:", big_matrix.right_kernel().basis(), flush=True)
\end{verbatim}
\normalsize
The output:
\tiny
\begin{verbatim}
Number of graphs in X: 210
Calculating X_formulas
Calculated X_formulas
Number of monomials in components of X: [101226, 101226, 101226, 101226]
Number of linearly independent graphs in X: 112
Maximal subset of linearly independent graphs: [0, 1, 2, 3, 4, 5, 6, 7, 8, 9, 10, 11, 13, 14, 15, 16, 17, 18, 19, 20, 21, 22, 23,
24, 25, 26, 27, 28, 30, 31, 32, 33, 34, 35, 36, 37, 38, 39, 41, 44, 45, 46, 47, 48, 49, 50, 51, 52, 53, 54, 55, 56, 57, 58, 59, 60,
61, 62, 63, 64, 65, 66, 67, 68, 70, 71, 72, 74, 75, 76, 78, 79, 80, 81, 84, 85, 86, 87, 89, 90, 91, 112, 113, 114, 124, 125, 137, 163,
164, 165, 166, 167, 168, 169, 170, 172, 173, 174, 175, 176, 180, 181, 183, 184, 185, 192, 193, 195, 200, 202, 203, 204]
Calculating adot
Calculated adot
Calculating X_a_formulas
Calculated X_a_formulas
Number of monomials in X_a_basis: 127344 127344
Calculating X_a_evaluation_matrix
Calculated X_a_evaluation_matrix
Calculating adot_vector
Calculated adot_vector
Calculating Q_tetra
Calculated Q_tetra
Calculating rhodot
Calculated rhodot
Have nice expression for Q_tetra: True
Calculating X_rho_formulas
Calculated X_rho_formulas
Number of monomials in X_rho_basis: 353532
Calculating X_rho_evaluation_matrix
Calculated X_rho_evaluation_matrix
Calculating rhodot_vector
Calculated rhodot_vector
Calculating X_solution_vector
Calculated X_solution_vector
X_solution_vector = (0, 12, 0, -12, 24, 0, -48, -24, -24, -24, 0, -24, 0, 0, -24, 0, 36, 0, 0, -18, 0, -24, 0, -24, -24, -48, -16, 0, 48,
0, 0, 0, 0, -16, -24, -16, -24, 12, 0, 0, 0, 24, -12, 0, 6, -16, -8, 0, 24, -16, 0, 0, 0, -24, 0, 0, 0, -48, -8, -12, 0, 0, -24, 0, 0, -24, 0, -24,
-24, 0, 0, 0, -8, -24, -16, -12, -48, -16, -24, -16, 0, 0, 0, 0, 0, 0, 0, -24, 0, 0, 0, -12, 0, 0, 0, -12, 0, 0, 0, -24,
-24, -12, 0, 0, 0, 0, 24, 0, 24, 0, 0, 0)
P.bracket(X_solution) == Q_tetra: True
Number of parameters in the solution: 7
Basis of kernel: [
(1, 0, 0, 0, 0, 1, -3, -2, 0, 0, 0, -1, 0, 0, 0, 0, 0, 1, 0, 0, 0, 0, 0, 1/2, 0, 0, 0, 0, 0, 0, 0, 0, 0, 0, -2, 0, -1, 2, -1/2,
0, 0, 0, 0, -1, 1, 0, 0, 2, -1, 0, -3, 0, 0, -1/2, 0, 0, 0, -3, 0, 0, 2, 0, -2, 2, 1, -2, 2, 0, -2, 0, 1/2, 2, 0, 0, 0, 0, -3, 0,
-2, 0, 2, 0, 1, -1, 0, -1/2, -1, 0, 0, 0, 0, 0, 0, 0, 0, 0, 0, 0, 0, 0, 0, 0, 0, 0, 0, 0, 0, 0, 0, 0, 0, 0),
(0, 1, 0, 0, 0, 0, -2, 0, 0, -2, 0, 0, 0, 0, 0, 0, 0, 0, 0, 0, 0, 0, 0, 0, 0, 0, 0, 0, 0, 0, 0, 0, 0, 0, 0, 0, 0,
0, 0, 0, 0, 0, 0, 0, 0, 0, 0, 0, 0, 0, 0, 0, 0, -1, 0, 0, 0, 0, 0, 0, 0, 0, 0, 0, 0, 0, 0, 0, 0, 0, 1, 0, 0, 0, 0, 0,
0, 0, 0, 0, 0, 0, 0, 0, 0, 0, 0, 0, 0, 0, 0, 0, 0, 0, 0, 0, 0, 0, 0, 0, 0, 0, 0, 0, 0, 0, 0, 0, 0, 0, 0, 0),
(0, 0, 0, 0, 1, 0, 0, 0, -1, 0, 0, 0, 0, 0, 0, 0, 0, 0, 0, 0, 0, 0, 0, 0, 0, 0, 0, 0, 0, 0, 0, 0, 0, 0, 0, 0, 0, 0, 0,
0, 0, 0, 0, 0, 0, 0, 0, 0, 0, 0, 0, 0, 0, 0, 0, 0, 0, 0, 0, 0, 0, 0, 0, 0, 0, 0, 0, 0, 0, 0, 0, 0, 0, 0, 0, 0, 0, 0, 0,
0, 0, 0, 0, 0, 0, 0, 0, 0, 0, 0, 0, -1/2, 0, 0, 0, -1/2, 0, 0, 0, 0, 0, -1/4, 1/2, 0, 0, 0, 0, 0, 1/2, 0, 0, 1/2),
(0, 0, 0, 0, 0, 0, 0, 0, 0, 0, 0, 0, 0, 0, 0, 0, 1, 0, 0, 0, 0, 0, 0, 0, 0, 0, 0, 0, 1, 0, 0, 0, 0, 0, 0, 0, 0, 0, 0, 0,
0, 1/2, 0, 1/2, 0, 0, 0, 0, 0, 0, 0, 0, 0, 0, 0, 0, 0, 0, 0, 0, 0, 0, 0, 0, 0, 0, 0, 0, 0, 0, 0, 0, 0, 0, 0, 0, 0, 0,
0, 0, 0, 0, 0, 0, 0, 0, 0, 0, 0, 0, 0, 0, 0, 0, 0, 0, 0, 0, 0, 0, 0, 0, 0, 0, 0, 0, 0, 0, 0, 0, 0, 0),
(0, 0, 0, 0, 0, 0, 0, 0, 0, 0, 0, 0, 0, 0, 0, 0, 0, 0, 0, 1, 0, 0, 0, 0, 0, 0, 0, 0, 0, 0, 0, 0, 0, 0, 0, 0, 0, 2, 0,
0, 0, 0, 0, 0, 1, 0, 0, 2, 0, 0, 0, 0, 0, 0, 0, 0, 0, 0, 0, 0, 0, 0, 0, 0, 0, 0, 0, 0, 0, 0, 0, 0, 0, 0, 0, 0, 0, 0, 0,
0, 0, 0, 0, 0, 0, 0, 0, 0, 0, 0, 0, 0, 0, 0, 0, 0, 0, 0, 0, 0, 0, 0, 0, 0, 0, 0, 0, 0, 0, 0, 0, 0),
(0, 0, 0, 0, 0, 0, 0, 0, 0, 0, 0, 0, 0, 0, 0, 0, 0, 0, 0, 0, 0, 1, 0, 1/2, 0, 0, 0, 0, 0, 0, 0, 0, 0, 0, 0, 0, 0, 0, -1/2,
0, 0, 0, 0, 0, 0, 0, 0, 0, -1, 0, -1, 0, 0, 0, 0, 0, 0, 0, 0, 0, 0, 0, 0, 0, 0, 0, 0, 0, 0, 0, 0, 0, 0, 0, 0, 0, 0, 0, 0,
0, 0, 0, 0, 0, 0, 0, 0, 0, 0, 0, 0, 0, 0, 0, 0, 0, 0, 0, 0, 0, 0, 0, 0, 0, 0, 0, 0, 0, 0, 0, 0, 0),
(0, 0, 0, 0, 0, 0, 0, 0, 0, 0, 0, 0, 0, 0, 0, 0, 0, 0, 0, 0, 0, 0, 0, 0, 0, 1, 0, 0, 0, 0, 0, 0, 0, 0, 0, 0, 0, 0, 0, 0,
0, 0, 0, 0, 0, 0, 0, 0, 0, 0, 0, 0, 0, 0, 0, 0, 0, 0, 0, 0, 1, 0, 0, 1, 0, 0, 1, 0, 0, 0, 0, 1, 0, 0, 0, 0, 0, 0, 0, 0, 1,
0, 0, 0, 0, 0, 0, 0, 0, 0, 0, 0, 0, 0, 0, 0, 0, 0, 0, 0, 0, 0, 0, 0, 0, 0, 0, 0, 0, 0, 0, 0)
]

###############################################################################
Hábrók Cluster

\end{verbatim}
\normalsize
\subsection{Expanding the 3D vanishing sunflower micro-graphs to 4D}\label{A8}

This result is condensed in Remark \ref{vanishingremark} and line 6 in Table~\ref{updowntable} on p.~\pageref{updowntable}.

The script:
\tiny
\begin{verbatim}
NPROCS=32

import warnings
warnings.filterwarnings("ignore", category=DeprecationWarning)

# Graphs for the vector field $X$

vanishing_expanded=[(0,2,4,7,1,3,5,8,4,2,6,9),(0,2,4,7,1,3,5,8,7,2,6,9),(0,5,4,7,4,3,5,8,1,2,6,9),(0,8,4,7,4,3,5,8,1,2,6,9),
(0,8,4,7,7,3,5,8,1,2,6,9),(0,5,4,7,7,3,5,8,1,2,6,9),(0,5,4,7,1,3,5,8,4,2,6,9),(0,8,4,7,1,3,5,8,4,2,6,9),(0,8,4,7,1,3,5,8,7,2,6,9),
(0,5,4,7,1,3,5,8,7,2,6,9),(0,2,4,7,1,3,5,8,4,5,6,9),(0,2,4,7,1,3,5,8,7,5,6,9),(0,2,4,7,1,3,5,8,7,8,6,9),(0,2,4,7,1,3,5,8,4,8,6,9)]

for i in [4,7]:
    for j in [5,8]:
        for k in [6,9]:
            vanishing_expanded.append((0,1,4,7,1,k,5,8,i,j,6,9))
            vanishing_expanded.append((0,2,4,7,1,k,5,8,i,j,6,9))
            vanishing_expanded.append((0,2,4,7,i,k,5,8,1,j,6,9))
            vanishing_expanded.append((0,j,4,7,i,k,5,8,1,2,6,9))
            vanishing_expanded.append((0,j,4,7,1,k,5,8,i,2,6,9))
            vanishing_expanded.append((0,1,4,7,i,k,5,8,1,j,6,9))
            
for i in [4,7]:
    for j1 in [5,8]:
        for j2 in [5,8]:
            vanishing_expanded.append((0,j1,4,7,1,3,5,8,i,j2,6,9))
            
for i1 in [4,7]:
    for i2 in [4,7]:
        for j1 in [5,8]:
            for j2 in [5,8]:
                vanishing_expanded.append((0,j1,4,7,i1,3,5,8,i2,j2,6,9))
                
for i1 in [4,7]:
    for i2 in [4,7]:
        for j1 in [5,8]:
            for j2 in [5,8]:
                for k in [6,9]:
                    vanishing_expanded.append((0,j1,4,7,i1,k,5,8,i2,j2,6,9))

X_graph_encodings = vanishing_expanded

# Convert the encodings to graphs:

from gcaops.graph.formality_graph import FormalityGraph
def encoding_to_graph(encoding):
    targets = [encoding[0:4], encoding[4:8], encoding[8:12]]
    edges = sum([[(k+1,v) for v in t] for (k,t) in enumerate(targets)], [])
    return FormalityGraph(1, 9, edges)

X_graphs = [encoding_to_graph(e) for e in X_graph_encodings]
print("Number of graphs in X:", len(X_graphs), flush=True)

# Formulas of graphs for $X$

# We want to evaluate each graph to a formula, inserting three copies of $\varrho\varepsilon^{ijk}$ 
# and three copies of the Casimir $a$ into the aerial vertices.

# First define the even coordinates $x,y,z$ (base), $\varrho, a$ (fibre) and odd coordinates $\xi_0, \xi_1, \xi_2$:

from gcaops.algebra.differential_polynomial_ring import DifferentialPolynomialRing
D4 = DifferentialPolynomialRing(QQ, ('rho','a1','a2'), ('x','y','z','w'), max_differential_orders=[3+1,1+3+1,1+3+1])
rho, a1, a2 = D4.fibre_variables()
x,y,z,w = D4.base_variables()
even_coords = [x,y,z,w]

from gcaops.algebra.superfunction_algebra import SuperfunctionAlgebra

S4.<xi0,xi1,xi2,xi3> = SuperfunctionAlgebra(D4, D4.base_variables())
xi = S4.gens()
odd_coords = xi
epsilon = xi[0]*xi[1]*xi[2]*xi[3] # Levi-Civita tensor

# Now evaluate each graph to a formula:

import itertools
from multiprocessing import Pool

def evaluate_graph(g):
    E = x*xi[0] + y*xi[1] + z*xi[2] + w*xi[3] # Euler vector field, to insert into ground vertex. Incoming derivative d/dx^i will result in xi[i].
    result = S4.zero()
    for index_choice in itertools.product(itertools.permutations(range(4)), repeat=3):
        sign = epsilon[index_choice[0]] * epsilon[index_choice[1]] * epsilon[index_choice[2]]
        # NOTE: This assumes the ground vertex is labeled 0, the vertices of out-degree 4 are labeled 1, 2, 3, 
        # and the Casimirs are labeled 4, 5, 6; 7, 8, 9.
        vertex_content = [E, S4(rho), S4(rho), S4(rho), S4(a1), S4(a1), S4(a1), S4(a2), S4(a2), S4(a2)]
        for ((source, target), index) in zip(g.edges(), sum(map(list, index_choice), [])):
            vertex_content[target] = vertex_content[target].derivative(even_coords[index])
        result += sign * prod(vertex_content)
    return result

print("Calculating X_formulas", flush=True)
X_formulas = []
with Pool(processes=NPROCS) as pool:
    X_formulas = list(pool.imap(evaluate_graph, X_graphs))
# NOTE: Necessary fixup of parents after multiprocessing.
X_formulas = [S4(X_formula) for X_formula in X_formulas]
print("Calculated X_formulas", flush=True)

print("How many vanishing micrographs?", X_formulas.count(0))

print("Which are the vanishing micrographs?", [k+1 for (k, X) in enumerate(X_formulas) if X == 0])
\end{verbatim}
\normalsize
The output:
\tiny
\begin{verbatim}
Number of graphs in X: 118
Calculating X_formulas
Calculated X_formulas
How many vanishing micrographs? 54
Which are the vanishing micrographs? [1, 2, 3, 4, 5, 6, 7, 8, 9, 10, 11, 12, 13, 14, 15, 16, 17, 18, 19, 20, 24, 30, 36, 42, 48, 54, 57, 58, 59, 60, 
61, 62, 63, 64, 65, 66, 67, 68, 69, 70, 71, 86, 87, 88, 97, 98, 99, 100, 105, 106, 107, 108, 117, 118]

###############################################################################
Hábrók Cluster

\end{verbatim}
\normalsize 

\subsection{The 5D sunflower descendant encodings and the computation of one 5D formula}\label{5D}

This result is contained in Table~\ref{table1} on p.~\pageref{table1}.

The script:
\tiny
\begin{verbatim}
NPROCS=32

import warnings
warnings.filterwarnings("ignore", category=DeprecationWarning)

# Make the 5D sunflower encodings

X_graph_encodings = []
# Frog components (with tadpole)
for i2 in [2,5,8,11]:
   for j1 in [1,4,7,10]:
       for j2 in [1,4,7,10]:
           for k in [3,6,9,12]:
               X_graph_encodings.append((0,1,4,7,10,j1,k,5,8,11,j2,i2,6,9,12))
# Ballerina components (without tadpole)
for i1 in [2,5,8,11]:
    for i2 in [2,5,8,11]:
        for j1 in [1,4,7,10]:
            for j2 in [1,4,7,10]:
                for k in [3,6,9,12]:
                    X_graph_encodings.append((0,i1,4,7,10,j1,k,5,8,11,j2,i2,6,9,12))

from gcaops.graph.formality_graph import FormalityGraph
def encoding_to_graph(encoding):
    targets = [encoding[0:5], encoding[5:10], encoding[10:15]]
    edges = sum([[(k+1,v) for v in t] for (k,t) in enumerate(targets)], [])
    return FormalityGraph(1, 12, edges)

X_graphs = [encoding_to_graph(e) for e in X_graph_encodings]

print("How many 5D sunflower descendants?", len(X_graphs))

from gcaops.algebra.differential_polynomial_ring import DifferentialPolynomialRing
D5 = DifferentialPolynomialRing(QQ, ('rho','a1','a2','a3'), ('x','y','z','w','v'), max_differential_orders=[3+1,1+4+1,1+4+1,1+4+1])
rho, a1, a2, a3 = D5.fibre_variables()
x,y,z,w, v = D5.base_variables()
even_coords = [x,y,z,w,v]

from gcaops.algebra.superfunction_algebra import SuperfunctionAlgebra

S5.<xi0,xi1,xi2,xi3,xi4> = SuperfunctionAlgebra(D5, D5.base_variables())
xi = S5.gens()
odd_coords = xi
epsilon = xi[0]*xi[1]*xi[2]*xi[3]*xi[4] # Levi-Civita tensor

import itertools
def evaluate_graph(g):
    E = x*xi[0] + y*xi[1] + z*xi[2] + w*xi[3] + v*xi[4] # Euler vector field, to insert into ground vertex. Incoming derivative d/dx^i will result in xi[i].
    result = S5.zero()
    for index_choice in itertools.product(itertools.permutations(range(5)), repeat=3):
        sign = epsilon[index_choice[0]] * epsilon[index_choice[1]] * epsilon[index_choice[2]]
        # NOTE: This assumes the ground vertex is labeled 0, the vertices of out-degree 4 are labeled 1, 2, 3, and the Casimirs are labeled 4, 5, 6; 7, 8, 9.
        vertex_content = [E, S5(rho), S5(rho), S5(rho), S5(a1), S5(a1), S5(a1), S5(a2), S5(a2), S5(a2), S5(a3), S5(a3), S5(a3)]
        for ((source, target), index) in zip(g.edges(), sum(map(list, index_choice), [])):
            vertex_content[target] = vertex_content[target].derivative(even_coords[index])
        result += sign * prod(vertex_content)
    return result

# We check how long it takes to compute 1 formula

formulas_1 = evaluate_graph(X_graphs[0])
\end{verbatim}
\normalsize
The output:
\tiny
\begin{verbatim}
How many 5D sunflower descendants? 1280

Used walltime                  :   08:23:59
\end{verbatim}

\end{normalsize}

\end{document}